\documentclass[11pt,A4]{article}

\usepackage{amsmath,amsfonts,amssymb,amsthm}
\usepackage{mathrsfs}
\usepackage{latexsym}
\usepackage[english]{babel}
\usepackage{graphicx}
\usepackage[usenames,dvipsnames]{xcolor}

\usepackage{xcolor}

\usepackage[normalem]{ulem} 

\newcommand{\red}{{}}

\newcommand{\blues}{\color{black} }

\newcommand{\teal}{ {} }

\renewenvironment{proof}[1][\proofname]{{\bfseries #1.} }{\qed}

\def\Cov{{\rm Cov\,}}

\newcommand{\field}[1]{\mathbb{#1}}

\newcommand{\gjmr}{\gamma_j^{(M,R)}}

\newcommand{\ljmr}{\lambda_j^{(M,R)}}

\newcommand{\SKMR}{\sqrt{ K^{(M,R)}} }
\newcommand{\SCMR}{\sqrt{ C^{(M,R)}} }
\newcommand{\CMR}{C^{(M,R)}}
\newcommand{\LBR}{L^2_{S(M)}(B_R)}
\newcommand{\KMR}{K^{(M,R)}}

\newcommand{\R}{\field{R}}

\newcommand{\e}{{\rm e}}

\newcommand{\WW}{{\mathbb W}}

\newcommand{\cC}{{\cal C}}

\newcommand{\var}{\operatorname{Var}}

\def\authors#1{{ \begin{center} #1 \vspace{0pt} \end{center} } \smallskip}
\def\institution#1{{\sl \begin{center} #1 \vspace{0pt} \end{center} } }

\def\title#1{{\huge\bf  \begin{center} #1 \vspace{0pt} \end{center}  } \smallskip}

\def\E{{\mathbb{ E}}}

\def\P{{\mathbb{P}}}

\def\blu{\textcolor{blue}}
\def\blu{{}}

\def\paref#1{(\ref{#1})}

\newtheorem{theorem}{Theorem}[section]
\newtheorem{proposition}[theorem]{Proposition}
\newtheorem{lemma}[theorem]{Lemma}

\newtheorem{defn}[theorem]{Definition}

\newtheorem{remark}[theorem]{Remark}

\begin{document}

\title{Small Scale CLTs for the Nodal Length of Monochromatic Waves}

\date{Mar 2020}

\authors{\large Gauthier Dierickx$^{(1)}$, Ivan Nourdin$^{(2)}$, \\ Giovanni Peccati$^{(2)}$ and Maurizia Rossi$^{(3)}$}

\institution{{\rm (1)} Fakult\"at f\"ur Mathematik, Ruhr-Universit\"at Bochum\\
{\rm (2)} Unit\'e de Recherche en Math\'ematiques, Universit\'e du Luxembourg \\  {\rm (3)} Dipartimento di Matematica e Applicazioni, Universit\`a di Milano-Bicocca}

\begin{abstract} We consider the nodal length $L(\lambda)$ of the restriction to a ball of radius $r_\lambda$ of a {\it Gaussian pullback monochromatic random wave} of parameter $\lambda>0$ associated with a Riemann surface $(\mathcal M,g)$ without conjugate points. Our main result is that, if  $r_\lambda$ grows slower than $(\log \lambda)^{\teal 1/25}$, then (as $\lambda\to \infty$) the length $L(\lambda)$ verifies a Central Limit Theorem with the same scaling as Berry's random wave model -- as established in Nourdin, Peccati and Rossi (2019). Taking advantage of some powerful extensions of an estimate by B\'erard (1986) due to Keeler (2019), our techniques are \blu{mainly} based on a novel intrinsic bound on the coupling of smooth Gaussian fields, that is of independent interest, \blu{and moreover allow us to improve some estimates for the nodal length asymptotic variance of pullback random waves in Canzani and Hanin (2016)}.
In order to demonstrate the flexibility of our approach, we also provide an application to phase transitions for the nodal length of arithmetic random waves on shrinking balls of the $2$-torus.

\medskip


\medskip


\end{abstract}


\section{Introduction}

Let $(\mathcal{M}, g)$ be a compact, \blu{smooth}, Riemannian surface \blu{without boundary}, and denote by $\phi_\lambda$ \blu{the} Gaussian {\bf monochromatic random wave} on $\mathcal{M}$ with parameter $\lambda>0$ ({see \S \ref{overview} for precise definitions}). The aim of the present paper is to study the local behaviour of the {\bf nodal set} of $\phi_\lambda$, when $\lambda\to \infty$ and $\phi_\lambda$ is restricted to a ball whose radius converges to zero \blu{as a function of $\lambda$}. 


\medskip

Our main result, stated in Theorem \ref{t:main} below, is that if $\mathcal{M}$ has \emph{no conjugate points} and $r_\lambda = o \left ( (\log \lambda )^{\teal 1/25}\right )$, then the {\bf nodal length} of the {\it pullback} wave $\phi_\lambda^{x_0}$ \blu{associated with $\phi_\lambda$ at a point $x_0\in \mathcal M$}
and restricted to a ball \blu{on the tangent space $T_{x_0}\mathcal M$} of radius $r_\lambda$, verifies a Central Limit Theorem (CLT) with exactly the same asymptotic behaviour of mean and variance as for {\bf Berry's random wave model} on $\R^2$ --- see \cite{berry77, berry2002, NPR}.

\medskip

Our techniques are based on three main tools: (i) a quantitative extension of an estimate by B\'erard \cite{Ber0} (see also \cite{bon}) \blu{due to Keeler \cite{Kee18}}, yielding an explicit bound on the rate of convergence of covariance functions of pullback waves \blu{on manifolds without conjugate points} (see Theorem \ref{t:newberard}), (ii) some new explicit estimates on the {\bf coupling of smooth Gaussian fields} in $C^k$ topologies (see Theorem \ref{t:wcoupling}), and (iii) an application of a {\bf mixed Kac--Rice formula} in order to control the discrepancy of nodal lengths associated with coupled random functions. 

\medskip

In particular, our way of exploiting the estimates at Point (ii) above will be to explicitly couple, in the $C^1$ topology, the pullback wave $\phi^{x_0}_\lambda$ with a copy of Berry's random field, in such a way that the CLT \blu{for the pullback nodal length} can be directly inferred from the main result of \cite{NPR}. {\blues We stress that the field $\phi^{x_0}_\lambda$ is, in general, {\it not} stationary: it follows that  -- in order to couple $\phi^{x_0}_\lambda$ with Berry's random waves -- we cannot take advantage of the techniques recently developed in \cite{BeMa}, that only apply to the coupling of stationary fields, via the optimal pairing of spectral measures with respect to quadratic transport distances. We also refer to Remark \ref{r:sodintrick} below for a brief comparison with a coupling technique outlined in \cite[Section 3.1.1]{S}}.

\medskip

It is also important to notice that our application of Kac--Rice at Point (iii) will allow us to deduce an exact asymptotic relation for the variance of the nodal length of $\phi^{x_0}_\lambda$. As we will see in more detail below, exact asymptotic characterisations for the mean and variance of nodal lengths (and, a fortiori, second order results like central and non-central limit theorems) are typically available only for exactly solvable models, like e.g. {\bf random spherical harmonics} \blu{\cite{Ber,Wig, MRW}}, {\bf arithmetic random wave}s \blu{\cite{RW,KKW, MPRW, DNPR, Cam, PR, BM17, BMW}}, or the already quoted Berry's planar waves \cite{berry2002, NPR}. To the best of our knowledge, our Theorem \ref{t:main} is the first exact second order result for nodal lengths of monochromatic waves holding for such a general class of random fields.

\medskip

As explained below, the results of the present paper provide a counterpart to the laws of large numbers proved by Canzani and Hanin in \cite[Theorem 1]{CH}, and also yield an explicit quantitative answer to a problem left open in \cite[Section 1.4.1]{NPR}. One should notice that the condition $r_\lambda = o\left ( (\log \lambda )^{\teal 1/25}\right )$ is much more restrictive than the requirement $r_\lambda = o(\lambda)$ that is \blu{sufficient} for the main results of \cite{CH} to hold: this is due to the fact that the rate $o\left ( (\log \lambda )^{\teal 1/25}\right )$ is the {\it maximal} one for which we can effectively couple the pullback wave of $\phi_\lambda$ with Berry's planar field in such a way that the difference between the corresponding normalized nodal lengths converge to zero in $L^2$ when $\lambda\to \infty$. While it is clear that the exponent ${\teal 1/25}$ is in part an artefact of some analytical inequalities that are applied in our proofs and could in principle be improved (see e.g. our use of Sobolev embedding in Section \ref{ss:coupling}), the logarithmic dependence on $\lambda$ is a consequence of \blu{\cite{Kee18} (see Theorem \ref{t:newberard})} refining a deep result by B\'erard \cite{Ber0}, and cannot easily be dispensed with -- see Section \ref{ss:maurizia}.

\medskip

In order to demonstrate the flexibility of our approach, we also provide an application to arithmetic random waves on the flat torus \blu{$\mathbb{T}^2:=\mathbb R^2/\mathbb Z^2$} -- which in principle do not enter the above described framework of pullback random waves -- yielding a small scale CLT that is {\blues related to} a conjecture of \blu{Benatar, Marinucci and Wigman, see \cite[\S 2.2]{BMW}}. In this case, we cannot directly apply the refinement of the estimate by B\'erard \cite{Ber0} (see also \cite{bon}) \blu{due to Keeler \cite{Kee18}}, and rely indeed on a direct argument based on a classical arithmetic estimate from \cite{KK} --- {\blues  the idea of using such an estimate for coupling arithmetic random waves with Berry's model already appears in \cite{BeMa}; see also \cite{Sartori}.

\medskip

We observe that similar problems for random spherical harmonics were studied by A. P. Todino in \cite{Tod18}. In such a paper, the author proves a CLT for the nodal length of random spherical harmonics in shrinking caps via a specific argument (a reduction principle). Our coupling techniques could be applied also in this framework, plausibly at the cost of worse estimates (in terms of conditions on the radius of the shrinking spherical cap) than those in \cite{Tod18}, because of the full generality of our approach.}

\medskip

We will now present a more detailed discussion of our main findings. In what follows, every random object is defined on a suitable probability space $(\Omega, \mathscr{F}, \P)$, with $\E$ and ${\bf Var}$ denoting respectively expectation and variance with respect to $\P$. Given two positive sequences $\{a_n\}, \{b_n\}$, we write $a_n\sim b_n$ if $a_n/b_n\to 1$.

\subsection{Overview and main results}\label{overview}

Let $(\mathcal{M}, g)$ be a compact, \blu{smooth}, Riemannian surface without boundary, and denote by $\Delta_g$ the associated Laplace-Beltrami operator. We write $\{f_j : j\geq 0\}$ to indicate an orthonormal basis of $L^2(\mathcal{M})$ composed of eigenfunctions of $\Delta_g$, that is,
$$
\Delta_g f_j +\lambda_j^2 f_j = 0, \quad j\geq 0,
$$
where the corresponding eigenvalues are such that $0=\lambda_0< \lambda_1\leq \lambda_2\leq\cdots\leq \lambda_j \nearrow \infty$. Following \cite{Z, CH}, we define the {\bf monochromatic random wave}  of parameter $\lambda>0$ on $\mathcal{M}$ to be the Gaussian random field on the manifold 
\begin{equation}\label{phi}
\phi_\lambda(x) := \frac{1}{\sqrt{\text{dim}(H_{\lambda})}} \sum_{\lambda_j\in [\lambda, \lambda + 1]} a_j f_j(x), \quad x\in \mathcal{M},
\end{equation}
where the $a_j$ are independent and identically distributed (i.i.d.) standard Gaussian random variables, and
$$
H_{\lambda} := \bigoplus_{\lambda_j \in [\lambda, \lambda + 1]} \text{Ker}(\Delta_g + \lambda_j^2\, \text{Id}),
$$
with the symbol $\text{Id}$ denoting the identity operator. The Gaussian field $\phi_\lambda$ is centred by construction, and its covariance kernel is given by 
\begin{equation}\label{covRiem}
K_{\lambda}(x,y) := \Cov(\phi_\lambda(x), \phi_\lambda(y)) = \frac{1}{\text{dim}(H_{\lambda})} \sum_{\lambda_j\in [\lambda, \lambda +1]} f_j(x) f_j(y),\quad x,y\in \mathcal M.
\end{equation}
``Short window'' random waves such as $\phi_\lambda$ in \eqref{phi} (for manifolds of arbitrary dimension) were first introduced by Zelditch \cite{Z} as general approximate models of random Gaussian Laplace eigenfunctions defined on manifolds not necessarily having spectral multiplicities; see e.g., \blu{\cite{CH, BW17, NS} and the refereces therein for further discussions}.

\medskip

Our aim in this paper is to study the local behaviour of the {\bf nodal set} of $\phi_\lambda$, as $\lambda \to \infty$, restricted to balls of decreasing radius. Our main tool in order to accomplish this task is the notion of a ``pullback'' random wave that we will study at {\bf points of isotropic scaling}.  In order to introduce these notions, we adopt the standard notation $J_0( r ) $, $r\geq 0$, to indicate the {\bf Bessel function of the first kind} with index 0, given by
$$
J_0(r) := \int_{S^1} e^{i \langle u, z \rangle} \, \frac{dz}{2\pi},
$$
where $\frac{dz}{2\pi}$ is the uniform \blu{probability} measure on the unit circle, and $u\in\R^2$ is any point such that $\|u\|=r$.

\medskip

Fix $x_0\in \mathcal M$, and consider the tangent space $T_{x_0}\mathcal M$ to the manifold at $x_0$: we define the {\bf pullback random wave} associated with $\phi_\lambda$ at $x_0$ as the Gaussian random field on $T_{x_0}\mathcal M$ given by 
$$
\phi_\lambda^{x_0}(u) := \phi_\lambda\left ( \exp_{x_0} \left ( \frac{u}{\lambda} \right )\right ),\qquad u\in T_{x_0} \mathcal M, 
$$
where $\exp_{x_0} : T_{x_0}\mathcal M \to \mathcal M$ is the exponential map at $x_0$. The planar field $\phi_\lambda^{x_0}$ is trivially centered and Gaussian and, by virtue of \paref{covRiem}, its covariance kernel $K^{x_0}_\lambda$ is given by
$$
K_{\lambda}^{x_0}(u,v) = K_{\lambda}\left(\exp_{x_0} \left ( \frac{u}{\lambda}  \right) , \exp_{x_0} \left ( \frac{v}{\lambda}\right ) \right ),\qquad u,v\in T_{x_0} \mathcal M.  
$$ 
A direct inspection of the above covariance kernel immediately shows that $\phi_\lambda^{x_0}$ is of class $C^\infty$ with probability one.

\begin{defn}[See \cite{CH}]{\rm We say that $x_0 \in \mathcal{M}$ is a point of {\bf isotropic scaling} if, for every positive function $\lambda \mapsto r_\lambda$ such that $r_\lambda = o(\lambda)$, as $\lambda\to \infty$, one has that
\begin{equation}\label{limit}
\sup_{u,v \in \mathbb{B}(r_\lambda)  } \Big| \partial^\alpha\partial^\beta \, \big\{  K_{\lambda}^{x_0}(u,v) - (2\pi) J_0(\| u-v\|_{g_{x_0}})\big\}\Big| \to 0,\quad \lambda\to \infty,
\end{equation}
where $\alpha, \beta\in \mathbb{N}^2$ are multi-indices labeling partial derivatives with respect to $u$ and $v$, respectively, $\| \cdot \|_{g_{x_0}}$ is the norm on $T_{x_0}\mathcal M$ induced by $g$, and $\mathbb{B}(r_\lambda)$ is the corresponding ball of radius $r_\lambda$ centred at the origin.
}
\end{defn}

\begin{remark}\label{r:berry} {\rm

\begin{enumerate}

\item[(a)] Sufficient conditions for a point $x_0$ to be of isotropic scaling are discussed e.g. in \cite[Section 2.5]{CH}, building on the findings \cite{CH0}. In particular, \cite[Theorem 1]{CH0} implies that a sufficient condition for $x_0\in \mathcal{M}$ to be of isotropic scaling is that the set
$$
\mathcal{L}_{x_0,x_0} := \{\xi \in S_{x_0} \mathcal{M} : \exists t>0 \mbox{ s.t. } \exp_{x_0}(t\xi) = x_0 \} 
$$
has volume 0 in $T_{x_0} \mathcal{M}$, where $S_{x_0} \mathcal{M}$ denotes the unit sphere in $T_{x_0} \mathcal{M}$ with respect to the norm $\|\cdot \|_{g_{x_0}}$. For every compact smooth manifold $\mathcal{M}$ and for every $x_0\in \mathcal{M}$, the property $| \mathcal{L}_{x_0,x_0} | = 0$ is generic in the space of all Riemaniann metrics \cite[Lemma 6.1]{SZ}. {It is also known that the condition $| \mathcal{L}_{x_0,x_0}| = 0 $ holds for every $x_0\in \mathcal{M}$ whenever $\mathcal{M}$ has no conjugate points} (and, in particular, when $\mathcal{M}$ is negatively curved).

\item[(b)] Relation \eqref{limit} implies that, for every $u,v\in T_{x_0}\mathcal{M}$ and every multi-indices $\alpha, \beta$, the two-dimensional Gaussian field $\{( \partial^\alpha\phi^{x_0}_\lambda (u), \partial^\beta\phi^{x_0}_\lambda (u)) : u\in T_{x_0}\mathcal{M}\}$ converges in the sense of finite-dimensional distributions to 
$$
\{\sqrt{2\pi}( \partial^\alpha\phi^{x_0}_\infty(u), \partial^\beta\phi^{x_0}_\infty (u)) : u\in T_{x_0}\mathcal{M}\},
$$   
where $\phi^{x_0}_\infty$ is the centered Gaussian field on $T_x\mathcal{M}$ with covariance 
$$
\E[\phi_\infty^{x_0} (u)\phi_\infty^{x_0} (v)] = J_0(\|u-v\|_{g_{x_0}}).
$$
One can easily check that, with probability one, $\phi_\infty^{x_0}$ is an eigenfunction with eigenvalue 1 of the Laplace operator on $T_{x_0}\mathcal{M}$ associated with the metric $g_{x_0}$. 

\item[(c)] ({\it Convention on the choice of coordinates}) Since in this paper we are only interested in second order results for a {\it fixed} $x_0\in \mathcal{M}$ of isotropic scaling, {we will always (tacitly) choose coordinates around $x_0$ in such a way that $g_{x_0} = {\rm Id}$}, and we will write $\|\cdot \|_{g_{x_0}} = \|\cdot \|$ in order to simplify the notation. In this way, the field $\phi_\infty^{x_0}$ at item (b) becomes universal (in the sense that it does not depend on $\mathcal{M}$) and can be identified with {\bf Berry's Random Wave Model} on $\R^2 \simeq T_{x_0}\mathcal{M}$. Such a field is defined as the unique (in distribution) centred real-valued random field $b = \{b(u ) : u \in \R^2\}$ such that $b$ is an eigenfunction of the Laplace operator $\Delta$ on $\R^2$ with eigenvalue 1, and $ b$ is {\it isotropic}, that is, the distribution of $b$ is invariant with respect to rigid motions of the plane. It can be proved that these requirements immediately imply that, necessarily, 
\begin{equation}\label{e:berryc}
\E[b(u)b(v)] = J_0(\|u-v\|);
\end{equation}
see \cite{NPR} for details. We observe that, when $g_x = {\rm Id}$, condition \eqref{limit} implies that, in the parlance of \cite{NS}, the ensemble $\{\phi^x_\lambda\}$ has {\it translation invariant local limits}.

\end{enumerate}
}
\end{remark}

As anticipated, the principal focus in our paper is the (random) {\bf nodal set} 
$$
(\phi_\lambda^{x_0})^{-1}(0)  := \{ u\in T_{x_0}\mathcal{M} : \phi_\lambda^{x_0}(u) = 0\}
$$
\blu{which is a.s.\ a smooth curve \cite{CH}}. 
For every $\lambda, r>0$ and $x_0\in \mathcal{M}$, we set
\begin{equation}\label{e:lprw}
L(\phi_\lambda^{x_0} ; r) := \mathcal{H}^1\left((\phi_\lambda^{x_0})^{-1}(0)\cap B_r\right),
\end{equation}
where $\mathcal{H}^1$ indicates the one-dimensional Hausdorff measure on $T_{x_0}\mathcal{M}\simeq \R^2$ and $B_r$ is the closed ball of radius $r$ centred at the origin; in other words, the random variable $L(\phi_\lambda^{x_0} ; r)$ represents the length of the restriction of the nodal set of $\phi_\lambda^{x_0}$ to $B_r$. Similarly, writing $b = \{b(u) : u\in \R^2\}$ for the Berry's random wave model defined in Remark \ref{r:berry}-(c) (in particular, formula \eqref{e:berryc}), we write
\begin{equation}\label{e:lbrw}
L(b ; r) : = \mathcal{H}^1\left(b^{-1}(0)\cap B_r\right).
\end{equation}
Note that ${\rm Vol}(B_r) = \pi r^2$. Our first statement is taken from \cite{CH}, and contains a `universal' law of large numbers for the nodal lengths $L(\phi_\lambda^x ; r)$ (observe that the convergence in $L^2$ at \eqref{e:lln} below is not stated in \cite[Theorem 1]{CH}, but it is rather an immediate consequence of the arguments in the proof).

\begin{theorem}[\bf Special case of Theorem 1 in  \cite{CH}]\label{t:ch} Let the above notation prevail and let $x_0$ be a point of isotropic scaling.
\begin{itemize}
\item[\rm (1)] For every fixed $r>0$, as $\lambda \to \infty$, one has that 
\begin{equation}\label{e:ptob}
L(\phi_\lambda^{x_0} ; r) \stackrel{law}{\longrightarrow} L(b ; r),
\end{equation}
where, here and for the rest of the paper, the symbol $\stackrel{law}{\longrightarrow}$ indicates convergence in distribution of random variables.

\item[\rm (2)] If the function $\lambda\mapsto r_\lambda$ is such that $r_\lambda = o(\lambda)$ as $\lambda\to \infty$, then
\begin{eqnarray}\label{e:lln}
\E\left[ \left( \frac{L(\phi_\lambda^{x_0} ; r_\lambda)}{r_\lambda^2} - \frac{\pi}{2\sqrt{2}}\right)^2\right] \rightarrow 0.
\end{eqnarray}
In particular, as $\lambda\to \infty$,
\begin{equation}\label{e:asmv}
\E\big[ L(\phi_\lambda^{x_0} ; r_\lambda)\big] \sim \frac{\pi}{2\sqrt{2} }r_\lambda^2, \,\,\,\mbox{and}\,\,\, {\bf Var}\big(L(\phi_\lambda^{x_0} ; r_\lambda)\big) = o(r_\lambda^4).
\end{equation}
\end{itemize}
\end{theorem}

In view of \eqref{e:lln}, the next logical step is to address the following question: {\em as $\lambda\to \infty$,  what is the nature of the fluctuations of $X_\lambda := L(\phi_\lambda^{x_0} ; r_\lambda)/{r_\lambda^2}$, around the limit $ \frac{\pi}{2\sqrt{2}}$? In particular, does a properly normalised version of $X_\lambda$ verify a \blu{CLT}?} 

\medskip

Plainly, answering such a question would require one to establish some non trivial lower bound for the function $$\lambda\mapsto {\bf Var}\big(L(\phi_\lambda^{x_0} ; r_\lambda)\big), \quad \lambda\to \infty,$$ and, in a generic setting like the one of Theorem \ref{t:ch}, such a task seems to be largely outside the scope of existing techniques -- see e.g. the discussion around Theorem 1 in \cite{CH}. 

\medskip

As anticipated, the main idea developed in the present paper is that, in the case of surfaces without conjugate points and if one considers mappings $\lambda \mapsto r_\lambda$ that diverge to infinity at a rate which is {\it considerably slower} than $\lambda$, then one can deduce precise informations about the fluctuations of $L(\phi_\lambda^{x_0} ; r_\lambda)$ from the following central limit theorem involving Berry's planar waves.

\begin{theorem}[\bf See \cite{berry2002} and \cite{NPR}] \label{t:npr} As $r\to \infty$, one has that
\begin{equation}\label{e:asmvb}
\E\big[ L(b ; r )\big] = \frac{\pi }{2\sqrt{2}}r^2, \,\,\,\mbox{and}\,\,\, {\bf Var}\big(L(b ; r)\big) \sim \frac{r^2\log r}{256}.
\end{equation}
Moreover, 
\begin{equation}\label{e:bclt}
\frac{ L(b ; r ) -  \E\big[ L(b ; r )\big]  }{ {\bf Var}\big(L(b ; r)\big)^{1/2}}\stackrel{law}{\longrightarrow} Z\sim \mathscr{N}(0,1),
\end{equation}
where $\mathscr{N}(0,1)$ indicates the one-dimensional Gaussian distribution with mean zero and variance 1.
\end{theorem}

The main achievement of our work is the following small scale second order result, establishing exact estimates for mean and variances, as well as a \blu{CLT}, for pullback random waves associated with manifolds having no conjugate points. As already recalled, this also answers a question left open in \cite[Section 1.4.1]{NPR}.

\begin{theorem}[\bf Small scale CLT for pullback random waves] \label{t:main} Let the above notation prevail, and assume that $(\mathcal{M}, g)$ is a compact, smooth, Riemannian surface without boundary and without conjugate points. Then, for every $x_0\in \mathcal{M}$ and every function $\lambda\mapsto r_\lambda$ such that $r_\lambda\to \infty$ and, as $\lambda\to\infty$, 
\begin{equation}\label{e:21}
\frac{ {\red r_\lambda^{ 25 } } }{(\log r_\lambda)^{\teal 4}} = o(\log \lambda)
\end{equation}
one has that
\begin{equation}\label{e:asmvpb}
\E\big[ L(\phi_\lambda^{x_0} ; r_\lambda) \big] \sim \frac{\pi\, r_\lambda^2}{2\sqrt{2}}, \quad \,\,\, {\bf Var}\big(L(\phi_\lambda^{x_0} ; r_\lambda)\big) \sim \frac{r_\lambda^2\log r_\lambda}{256},
\end{equation}
and  
\begin{equation}\label{e:pbclt}
\frac{ L(\phi_\lambda^{x_0} ; r_\lambda) -  \E\big[ L(\phi_\lambda^{x_0} ; r_\lambda)\big]  }{ {\bf Var}\big(L(\phi_\lambda^{x_0} ; r_\lambda)\big)^{1/2}}\stackrel{law}{\longrightarrow} Z\sim \mathscr{N}(0,1). 
\end{equation}

\end{theorem}

\begin{remark}\label{r:discussion}{\rm 

\begin{itemize}

\item[(a)] In the regime \eqref{e:21}, the second relation in \eqref{e:asmvpb} largely improves the estimate ${\bf Var}\big(L(\phi_\lambda^{x_0} ; r_\lambda)\big) = o(r_\lambda^4)$, which is valid for generic $r_\lambda = o(\lambda)$ -- see Theorem \ref{t:ch}. Note that \eqref{e:21} is implied by 
$r_\lambda = o ((\log \lambda)^{ {\red 1/{25 } } })$.

\item[(b)] Our proof of Theorem \ref{t:main} will be achieved along the following route: {\bf (i)} \blu{taking advantage of the new bound for the rate of convergence \blu{to zero} in \eqref{limit} in the case of manifolds without conjugate points \cite{Kee18} (see Theorem \ref{t:newberard})}, {\bf (ii)} using the quantitative information at Point {\bf (i)} in order to build a coupling of $\phi_\lambda^{x_0}$ and Berry's planar wave $b$ in such a way that, as $\lambda\to\infty$, the difference  $  \phi_\lambda^{x_0} - \sqrt{2\pi}\, b $ converges to zero (say, in $L^1(\P)$) in the $C^1$ topology of $B_{r_\lambda}$, when $r_\lambda \to \infty$ sufficiently slow (see Theorem \ref{t:wcoupling}), and {\bf (iii)} applying a `mixed Kac--Rice formula' in order to show that, for 
$r_\lambda = o( (\log \lambda)^{ {\red 1/25 } })$ 
and for the coupling of $\phi_\lambda^{x_0}$ and  $b$ at Point {\bf (ii)} one has actually that, as $\lambda\to\infty$,
$$
\left| \frac{ L(b ; r_\lambda ) -  \E\big[ L(b ; r_\lambda )\big]  }{ {\bf Var}\big(L(b ; r_\lambda)\big)^{1/2}} - \frac{ L(\phi_\lambda^{x_0} ; r_\lambda) -  \E\big[ L(\phi_\lambda^{x_0} ; r_\lambda)\big]  }{  {\bf Var}\big(L(b ; r_\lambda)\big)^{1/2}}\right| \longrightarrow 0 \,\, \mbox{ in } L^2(\P);
$$
see Section \ref{s:kr}.

\item[(c)] Assume that, for some coupling of $\phi_\lambda^{x_0}$ and $b$ one has that $\phi_\lambda^{x_0} -\sqrt{2\pi}\, b$ converges to zero in probability, with respect to the $C^1$ topology of $B_{r_\lambda}$, with $r_\lambda \to \infty$, that is: for every $\epsilon>0$,
$$
\mathbb{P}\left[ \max_{\alpha : |\alpha|\leq 1} \sup_{z\in B_{r_\lambda}} \big| \partial^\alpha \phi_\lambda^{x_0}(z) - \partial^\alpha \sqrt{2\pi}\, b \big| > \epsilon \right]  \rightarrow 0, \quad \lambda\to\infty.
$$
Then, in general, it is {\it not} possible to conclude that the difference $L(b ; r_\lambda ) - L(\phi_\lambda^{x_0} ; r_\lambda)$ also converges to zero in probability \blu{({this is in contrast with the case of a fixed radius ball} -- see e.g. \cite{APP})}. This observation explains the necessity of Step {\bf (iii)} in the strategy outlined at the previous item.

\end{itemize}

}
\end{remark}

\begin{remark}\rm 
{It is worth noting that our `coupling' approach to limit theorems and variance estimates can be in principle applied to the case of parametric Gaussian ensembles on manifolds, that are locally converging to translation invariant Gaussian fields -- see e.g. the general framework outlined in \cite[Section 1.2]{NS} -- as soon as limit theorems for nodal lengths (or more general functionals) associated with the latter are known. However, deducing variance asymptotics and limit theorems for local geometric functionals of generic stationary fields, similar to Theorem \ref{t:npr}, would require a {remarkable amount of technical work and novel ideas, and will be investigated elsewhere}. 
}
\end{remark}

\subsection{The case of arithmetic random waves}

The framework of pullback random waves described in the previous section does not encompass the case of some exactly solvable models of Gaussian Laplace eigenfunctions defined on manifolds having spectral multiplicities, such as the model of arithmetic random waves  \cite{RW, KKW, MPRW, PR} and random spherical harmonics \cite{Ber, Wig, MRW}. The techniques developed in this paper can nonetheless be suitably extended in order to deal with specific models of this type. The aim of this subsection (and of Section \ref{s:proofarw} below, containing the proof of our main Theorem \ref{t:arw}) is to state and prove a small scale CLT for arithmetic random waves restricted to fast shrinking ball. As discussed below, such a theorem represents a counterpart to the main findings in \cite{BMW}, and corroborates a conjecture stated therein \blu{\cite[\S 2.2] {BMW}}.

\subsubsection{Definitions and reminders on global results} 

It is well-known that the eigenvalues of the Laplace operator on the flat 2-torus $\mathbb T^2$ are of the form $-E_n$, where $E_n := 4\pi^2n$ and
$$
n\in S:= \{n\in \mathbb{Z} : n = a^2+b^2, \, a,b\in \mathbb{Z}\}
$$
is the set of integers that can be represented as the sum of two squares. For $n\in S$, denote by $\Lambda_n$ the set of frequencies 
$$
\Lambda_n = \lbrace \xi\in \mathbb Z^2 : \| \xi\|= \sqrt n\rbrace 
$$
and by $\mathcal N_n$ the cardinality of $\Lambda_n$ (that is, $\mathcal N_n$ is the multiplicity of the eigenspace corresponding to $-E_n$). For $n\in S$, consider the probability measure $\mu_n$ induced by $\Lambda_n$ on the unit circle $\mathbb S^1$:
$$
\mu_n = \frac{1}{\mathcal N_n} \sum_{\xi\in \Lambda_n} \delta_{\xi/\sqrt n}.
$$ 
Following \cite{RW}, for $n\in S$, the toral random eigenfunction $T_n$ (or {\bf arithmetic random wave} of order $n$) is defined as the centered Gaussian field on the torus with the following covariance function: for $x,y\in \mathbb T^2$,
\begin{equation}\label{cov_int}
\Cov(T_n(x), T_n(y)) = \frac{1}{\mathcal N_n} \sum_{\xi\in \Lambda_n} \e^{i 2\pi\langle \xi, x-y\rangle } = \int_{\mathbb S^1} \e^{i2\pi \sqrt{n} \langle \theta, x-y\rangle}\,d\mu_n(\theta).
\end{equation}
It is easily checked that, with probability one, $\Delta T_n = -4\pi^2n T_n$, that is, $T_n$ is an eigenfunction of $\Delta$ with eigenvalue $-E_n$. As discussed in $\cite{KKW} $, there exists a density-$1$ subsequence $\lbrace n_j : j\geq 1\rbrace\in S$ such that, as $j\to +\infty$, 
$$
\mu_{n_j} \Rightarrow \frac{dz}{2\pi},
$$
where $\frac{dz}{2\pi}$ denotes as before the uniform \blu{probability} measure on the unit circle, and $\Rightarrow$ stands for weak convergence. \blu{For this subsequence, for $x,y\in \mathbb T^2$, 
$$
\Cov\left (T_{n_j}(x/2\pi\sqrt{n_j}\right ), T_{n_j}(y/2\pi\sqrt{n_j})) \to  \int_{\mathbb S^1} \e^{i \langle z , x-y\rangle}\,\frac{dz}{2\pi} = J_0(\|x-y\|),
$$
i.e. the scaling limit of $T_{n_j}$ is Berry's RWM.}

Let us now set $\mathcal L_n := \text{length}(T_n^{-1}(0))$. The expected nodal length was computed in \cite{RW} to be equal to
$$
\E[\mathcal L_n] = \frac{1}{2\sqrt 2} E_n,
$$
while in \cite{KKW} it is shown that, as $\mathcal N_n\to +\infty$, the variance of $\mathcal L_n $ satisfies the following exact relation
$$
\var(\mathcal L_n) \sim \frac{1+\widehat{\mu_n}(4)^2}{512} \, \frac{E_n}{\mathcal N_n^2},
$$
where $\widehat{\mu_n}(4)$ denotes the fourth Fourier coefficients of $\mu_n$. In order to have an asymptotic law for the variance, one should select a subsequence $\lbrace n_j\rbrace$ of energy levels such that (i) $\mathcal N_{n_j}\to +\infty$ and (ii) $|\widehat{\mu_n}(4)|\to \eta$,  for some $\eta\in [0,1]$. Note that for each $\eta\in [0,1]$, there exists a subsequence $\lbrace n_j\rbrace$ such that both (i) and (ii) hold (see \cite{KKW, KW}). For such subsequences, the asymptotic distribution of the nodal length was shown to be non-Gaussian in \cite{MPRW}: 
\begin{equation}\label{NCLT}
\frac{\mathcal L_{n_j} - \E[\mathcal L_{n_j}]}{\sqrt{\var(\mathcal L_{n_j})}}\stackrel{law}{\longrightarrow} \frac{1}{2\sqrt{1+\eta^2}} (2 - (1-\eta) Z_1^2 - (1+\eta)Z_2^2),
\end{equation}
where $Z_1$ and $Z_2$ are i.i.d. standard Gaussian random variables. A complete quantitative version (in Wasserstein distance) of \paref{NCLT} is given in \cite{PR}.  

\subsubsection{Phase transitions for nodal lengths on shrinking balls}

The following remarkable statement is extrapolated from \cite[Theorem 1.1]{BMW} \blu{and} contains a characterization of the fluctuation of the nodal length of arithmetic random waves above the Planck scale. For every $n\in S$ and every $r>0$, we set
$$
\mathcal{L}_n(r):=  \text{length}\Big(T_n^{-1}(0) \cap B_{r /\sqrt{n} }\Big),
$$
where $B_r$ denotes the ball of radius $r$ centred at the origin. Note that, in general
$$
\E\left[ \text{length}\Big(T_n^{-1}(0) \cap B_{r }\Big) \right] = \frac{\pi r^2}{2\sqrt{2}} E_n.
$$

\begin{theorem}[\bf Special Case of Theorem 1.1 in \cite{BMW}]\label{t:bmwarw} For every $\gamma\in (0,1/2) $, there exists a density one sequence $\{n_j\}\subset S$ verifying the following properties:
\begin{enumerate}
\item[\rm 1.] as $n_j\to \infty$, one has that $\mathcal{N}_{n_j}\to \infty$ and $\mu_{n_j}$ converges weakly to the uniform measure on $\mathbb{S}^1$;

\item[\rm 2.] as $n_j\to \infty$,
$$
{\bf Var}( \mathcal{L}_{n_j} (n_j^\gamma)) \sim \frac{1+\widehat{\mu_{n_j} }(4)^2}{512} \, \frac{E_{n_j}}{\mathcal N_{n_j}^2}\times  \big\{ \pi (n_j^{\gamma-1/2})^2 \big \}^2 ;
$$

\item[\rm 3.] as $n_j\to \infty$,
$$
\frac{\mathcal L_{n_j} - \E[\mathcal L_{n_j}]}{\sqrt{\var(\mathcal L_{n_j})}}\stackrel{law}{\longrightarrow} \blu{1 -\frac{Z_1^2 + Z_2^2}{2}},
$$
where $Z_1,Z_2$ are two independent standard Gaussian random variables.

\end{enumerate}

\end{theorem}

In \cite[Section 2.2]{BMW} a general conjecture is stated concerning nodal lengths of arithmetic random waves, containing in particular the following

\smallskip

\noindent{\bf Conjecture.} {\it There exists $A_0>0$ such that {\rm (i)} the conclusion of Theorem \ref{t:bmwarw} continues to hold if one replaces the sequence $n^\gamma$ with any sequence $\alpha_n \geq (\log n )^C$, for any $C>A_0$, and {\rm (ii)} the conclusion of Theorem \ref{t:bmwarw} fails to hold if one replaces the sequence $n^\gamma$ with any sequence $\alpha_n = (\log n )^C$, for any $C< A_0$}.

\medskip 

The following statement is the main result of this section and shows that, if such an $A_0$ exists, then necessarily $A_0\geq \frac{1}{18 }(\log\pi - \log 2) = 0.02508 ...$ .

\begin{theorem}\label{t:arw}Fix $\rho <\frac{1}{2}(\log\pi - \log 2) = 0.225791...$. Then there exists a density one sequence $\{n_j\}\subset S$ such that, as $n_j\to \infty$,
\begin{enumerate}
\item[\rm 1.] $\mathcal{N}_{n_j}\to \infty$ and $\mu_{n_j}$ converges weakly to the uniform measure on $\mathbb{S}^1$;

\item[\rm 2.] for every sequence $n\mapsto\alpha_n$ such that $\alpha_n = O( (\log n)^{\rho/9})$,
$$
{\bf Var}( \mathcal{L}_{n_j} (\alpha_{n_j} )) \sim \frac{1}{n_j} \frac{\alpha_{n_j}^2\log \alpha_{n_j}}{256}
$$
and
$$
\frac{\mathcal L_{n_j}(\alpha_{n_j}) - \E[\mathcal L_{n_j}(\alpha_{n_j})]}{\sqrt{\var(\mathcal L_{n_j}(\alpha_{n_j}))}}\stackrel{law}{\longrightarrow} Z,
$$
where $Z$ denotes as before a standard Gaussian random variable.
\end{enumerate}

\end{theorem}

{\blues

\begin{remark}{\rm

\begin{itemize}

\item[(a)] The conclusion of Theorem \ref{t:arw} is implicitly based on a more general estimate, yielding that if $L(n)$ denotes the nodal length of the rescaled random wave $x\mapsto T_n (x/2\pi\sqrt{n})$ on the ball with radius $\alpha_n$ and $L'(n)$ denotes the nodal length of Berry's random wave on the same ball, then one can couple each $L(n_j)$ and $L'(n_j)$ in such a way that
$$
\E[ (L(n_j) - L'(n_j))^2] = O\left( \frac{\alpha_{n_j}^5}{  ( \log n_j )^{\rho/3}}\right).
$$

\item[(b)] While circulating an earlier draft of the present paper, it was brought to our attention that comparable lower bounds on the constant $A_0$ have been independently obtained in \cite[Theorem 1.4]{Sartori}, by combining coupling techniques from \cite{BeMa} with explicit estimates of the nodal length of perturbed random fields.

\end{itemize}

}

\end{remark}
}


\subsection{Plan}

{In \S \ref{Sestimates} we first recall the new result by Keeler on wave equation theory on compact manifolds without conjugate points, which improves some estimates by B\'erard, and then we state our result on coupling of Gaussian fields. In \S \ref{s:kr} we prove our main theorem, dealing first with an application of a mixed Kac--Rice formula in order to control the discrepancy of nodal lengths associated with coupled random functions. Some technical lemmas are collected in \S \ref{Stechnical}. Finally in \S \ref{s:proofarw}, we prove our main result on the phase transition for nodal lengths of arithmetic random waves.}

\section{Some estimates}\label{Sestimates}

\subsection{Explicit rates of convergence on manifolds without conjugate points}\label{ss:maurizia}

As written in Remark \ref{r:berry}, if $(\mathcal M, g)$ has no conjugate points, then every $x_0\in \mathcal M$ is of isotropic scaling and \paref{limit} holds at any point.  However in order to prove our main result, we need an explicit rate in \paref{limit}. \\
In \cite{Ber0}, the author solves this question for the on-diagonal case, i.e., for $u=v$ (using the same notations as in \paref{limit}) and when no derivatives are involved. The recent result by Keeler in  \cite{Kee18} is a breakthrough in this direction, indeed he greatly improved the error in Weyl's law on manifolds without conjugate points considering also the case of off-diagonal terms and derivatives of all order. \\
The following theorem is Corollary 1.1 in \cite{Kee18} and completely answers the question addressed just above, giving a (logarithmic) rate for \paref{limit} in full generality.
\begin{theorem}[Corollary 1.1 in \cite{Kee18}]\label{t:newberard}
Let $(\mathcal M, g)$ be a smooth, compact, Riemannian manifold of dimension two without conjugate points, then as $\lambda \to +\infty$, for any multi-indices $\alpha,\beta\in \mathbb N^2$
\begin{equation*}
\sup_{u,v \in \mathbb{B}(r_\lambda)  } \Big| \partial^\alpha\partial^\beta \, \big\{  K_{\lambda}^{x_0}(u,v) - (2\pi) J_0(\| u-v\|)\big\}\Big|  = O\left(\frac{1}{\log \lambda}\right),
\end{equation*}
whenever $r_\lambda = O\left (\sqrt{\frac{\lambda}{\log \lambda}} \right )$. Here the implicit constant in the $O$-notation depends on the choice of $x_0\in \mathcal M$ and $r_\lambda$, and on the order of differentiation. 
\end{theorem}
Note that $\mathbb{B}(r_\lambda)$ corresponds to a shrinking ball of radius $\frac{r_\lambda}{\lambda} = O\left ( \frac{1}{\sqrt{\lambda \log \lambda}}  \right )$ on $\mathcal M$.

\subsection{Coupling of smooth Gaussian fields}\label{ss:statcoupling}

We now state our main results about the coupling of smooth Gaussian fields on subsets of $\R^d$. Since the present paper only involves infinitely differentiable random fields, we will uniquely focus on the case of covariance functions of class $\cC^{\infty, \infty}$; it is a standard task to adapt our findings to the case of covariance functions of class $\cC^{k,k}$, for some finite integer $k$. Proofs are deferred to Section \ref{ss:hsgrow} and Section \ref{ss:coupling}.

\medskip

For the rest of the section, fix an integer $d\geq 1$. For every $R\geq 1$, we denote as before by $B_R$ the open ball centered at the origin and with radius $R$, and write $|B_R|$ for the volume of $B_R$. 
For integers $p,q\geq 1$, we denote by $\WW^{p,q}(B_R)$ and $\cC_b^p(B_R)$, respectively, the Sobolev space of $B_R$ with indices $p,q$, and the Banach space of continuous functions on $B_R$ having partial derivatives of order $\leq p$ that are uniformly continuous on $B_R$.

\medskip

In what follows, we shall consider two real-valued covariance kernels
$$
K : (x,y) \mapsto K(x,y), \quad \mbox{and} \quad C : (x,y)\mapsto C(x,y)
$$
defined on $\R^d \times \R^d$; we assume that $C$ is of class $\cC^{\infty, \infty}$ ($C$ is infinitely continuously differentiable in each variable $x$ and $y$), and $K$ is of class $\cC_b^{\infty, \infty}$ ($K$ is infinitely continuously differentiable in each variable $x$ and $y$, with bounded derivatives of every order). Standard results (see e.g. \cite[Section 1.4]{AT} or \cite[Appendix A.9]{NS}) imply that, if $X$ is a centered Gaussian field on $\R^d$ with covariance give by $K$ or $C$, then $X$ admits a modification that is of class $\cC^{\infty}$ with probability one; in what follows, we will uniquely (tacitly) consider such a modification.

\medskip

Given multi-indices $\alpha, \beta$, we introduce the shorthand notation
$$
K_{\alpha\beta} (x,y) := \partial^\alpha\partial^\beta K(x,y), \quad x,y\in\R^d, 
$$
and we define analogously the kernel $C_{\alpha\beta}$. For every integer $M=0,1,2,...$, we will write
\begin{equation}\label{e:sm}
S(M) : = \{\alpha : \alpha \mbox{ is a multi-index s.t. } |\alpha|\leq M\}.
\end{equation}
The following statement contains a coupling result for Gaussian fields belonging to Sobolev spaces, and is one of the main tools exploited in the sections to follow.

\begin{theorem}\label{t:wcoupling} Let the above notation and assumptions prevail, fix $M=0,1,...$ and $R\geq 1$, and write
\begin{equation}\label{e:eta}
\eta = \eta(M,R) := \max_{\alpha, \beta\in S(M)} \sup_{x,y\in B_R} \Big| K_{\alpha\beta}(x,y) - C_{\alpha\beta}(x,y)\Big|.
\end{equation}
Then, on some probability space $(\Omega_0, \mathscr{F}_0, \P_0)$, there exists a centred two-dimensional Gaussian field
$$
\{(X_0(z),Y_0(z) ) : z\in B_R\} 
$$
such that $X_0$ has covariance $K$, $Y_0$ has covariance $C$ and
\begin{equation}\label{e:sobcoupling}
\E_0\left[\| X_0-Y_0\|^2_{\WW^{M,2}(B_R)} \right]\leq A \Big\{\eta\,  |B_R| + \sqrt{\eta}\,  |B_R|^{1/2} R^{\frac{3d+1}{2}}\Big \},
\end{equation}
where $A = A(M,d)$ is an absolute finite constant independent of $R$. If moreover $M> j:= \lfloor \frac{d}{2}\rfloor +1$, then one has the additional estimate
\begin{equation}\label{e:bancoupling}
\E_0\left[\| X_0-Y_0\|^2_{\cC^{M-j}_b(B_R)} \right]\leq A' \, R^{2M-d} \Big\{\eta\,  |B_R| +\sqrt{\eta}\,  |B_R|^{1/2} R^{\frac{3d+1}{2}}\Big \},
\end{equation}
where $A' = A'(M,d)\in (0, \infty)$ is independent of $R$.
\end{theorem}

The estimate \eqref{e:bancoupling} is deduced by combining \eqref{e:sobcoupling} with a version of the Sobolev embedding theorem for open sets of $\R^d$ such as the one stated in \cite[Theorem 2.7.2]{DD} --- see Section \ref{ss:coupling} for a detailed proof.

{\blues

\begin{remark}\label{r:sodintrick}{\rm There is an alternate procedure for coupling smooth Gaussian processes via a discretization procedure, as described in \cite[Section 3.1.1]{S}. Such a procedure consists in the following steps: (i) for some parameter $\alpha>0$, choose an $\alpha$-net contained in $B_R$, (ii) build an optimal coupling of $\phi_\lambda^x$ and $\sqrt{2\pi}\, b$ on the $\alpha$-net fixed above (using some optimal criterion for coupling of Gaussian vectors -- see e.g. \cite{OP}), (iii) extend the finite coupling at Step (ii) by using an additional collection of independent Gaussian random variables, (iv) compute a bound on the $\cC^1_b(B_R)$ distance between the coupled fields by using some a priori estimates on their $\cC^2$ norms, combined e.g. with some version of the Kolmogorov-Landau inequality on a finite domain (such as e.g. an appropriate tensorization of \cite[Theorem 3.5]{chen}), and optimize in $\alpha$. While preparing our work, we actually pursued such a strategy in full detail, and managed to obtain a bound analogous to the content of Theorem \ref{t:wcoupling}, but where the right-hand side of \eqref{e:bancoupling} is replaced by the quantity 
\begin{equation}\label{e:boundzz}
A \sqrt{ \eta^{1/(d+1)}   R^{(3d+1)/(d+1)} (\log R)^{(2d+1)/(d+1)}}.
\end{equation}
Using such a bound in our proof yields a version of Theorem \ref{t:main} where condition \eqref{e:21} is replaced by the slightly stronger requirement that 
\begin{equation}\label{e:boundz}
\frac{ {\red r_\lambda^{ 28 } } }{(\log r_\lambda)^{\teal 7}} = o(\log \lambda)
\end{equation}
We also mention that, with respect  to the methods developed in the present paper, the approach of \cite{S} has the advantage of allowing one to deal directly with random fields of class $C^2$. It is also reasonable to expect that \eqref{e:boundzz} might perform better than \eqref{e:bancoupling} for large values of the dimensional parameter $d$. }

\end{remark}

}

\section{Proof of Theorem \ref{t:main}}\label{s:kr}

\subsection{Preparation}
Fix $x\in \mathcal{M}$. For the rest of the Section we write $K(x,y) = K(x-y)= (2\pi) J_0(\|x-y\|)$ and $C_\lambda(x,y) = K_\lambda^x (x,y)$. Fix $r_\lambda = o\left ( (\log \lambda)^{1/25 }\right )$. 
By virtue of Theorem \ref{t:newberard}, we know that
\begin{equation}\label{e:berardino}
\eta_\lambda := \max_{\alpha,\beta\in S(3)} \sup_{x,y\in B_{r_\lambda}} \left| \partial^\alpha \partial^\beta \{ K(x - y) - C_\lambda(x,y) \} \right| = O\left(\frac{1}{\log \lambda}\right),
\end{equation}
where the notation $ \partial^\alpha \partial^\beta K(x - y)$ indicates that the operator $\partial^\alpha$ acts on the variable $x$, and $\partial^\beta$ on the variable $y$. According to Theorem \ref{t:wcoupling}, as applied to the case $d=2$ and $M=3$, for every $\lambda>0$ there exists a jointly Gaussian coupling $(Y_\lambda, X)$ of $\phi_\lambda^x$ and $\sqrt{2\pi}\,  b$ such that
\begin{equation}\label{def_a_lambda}
\E\left[\| X-Y_\lambda\|^2_{\cC^1_b(B_{r_\lambda})}\right] \leq A\, \sqrt{\frac{r_\lambda^{17}   }{\log \lambda}} =: a(\lambda)\to 0, 
\end{equation}
where the constant $A$ is independent of $\lambda$ (we stress that the probability space $(\Omega_\lambda, \mathscr{F}_\lambda, \P_\lambda)$ on which the coupling is defined depends in general on $\lambda$, but we will omit such a dependence for the sake of readability). For every $\lambda>0$ and every $x,y\in B_R$ we introduce the following notation for mixed covariances: for every multi-indices $\alpha, \beta$
\begin{equation} \label{Eq_MixedCovKernel}
M^\lambda_{\alpha,\beta}(x,y) := \E[\partial^\alpha X(x) \partial^\beta Y_\lambda(y)]
\end{equation}
{\red and 
\begin{equation}\label{eta}
\zeta_{ \lambda} := {\teal \max_{\alpha,\beta\in S(3)}} \sup_{x,y\in B_{r_\lambda}} \left| \partial^\alpha \partial^\beta K(x - y) - M^\lambda_{\alpha,\beta}(x,y)  \right|.
\end{equation}
The previous discussion implies that, for every $\alpha, \beta \in S(3)$, there exists a finite constant $B$, independent of $\lambda$, such that, for every $\alpha, \beta \in S(3)$ and every $x,y\in B_{r_\lambda}$,
\begin{equation}\label{e:mixcov}
\left| M^\lambda_{\alpha,\beta}(x,y) - \partial^\alpha\partial^\beta K(x,y)\right| \leq \E[|\partial^\alpha X(x)|\cdot |\partial^\beta Y_\lambda(y) - \partial^\beta X(y)|]\leq B\sqrt{a(\lambda)},
\end{equation}
{where $a(\lambda)$ is defined in \paref{def_a_lambda}}.
Now we adopt a strategy close to the one {pursued} in \cite[Section 7.1]{NPR}, which is in turn inspired by \cite{ORW, RW}. We fix a large parameter $N>0$ (independent of $\lambda$, and whose value will be clarified later). We denote by $Q_0$ the square $[0, 1/N)^2$ and denote by $Q_{\bf z}$ the translation of $Q_0$ in the direction ${\bf z}/N$, where ${\bf z}\in \mathbb{Z}^2$. We write $\mathcal{Q}$ for the collection of all $Q_{\bf z}$ and, for every $\lambda>0$ we set $\mathcal{Q}_\lambda := \{Q_{\bf z} : Q_{\bf z} \cap B_{r_\lambda}\neq \emptyset\}$, in such a way that $| \mathcal{Q}_\lambda | = O(r_\lambda^2)$, as $\lambda\to \infty$, where the constant implicitly involved in such a relation only depends on the choice of $N$.
Fix a small number $\epsilon>0$. 
\begin{defn}
We say that two cubes $Q_{\bf x}$ and $Q_{\bf y}$ are {\it singular} if there exists $(x,y)\in Q_{\bf x}\times Q_{\bf y}$ such that, for some $\alpha, \beta\in S(1)$, $K_{\alpha\beta}(x,y) >\epsilon$. 
\end{defn}
We will need the following technical lemma. 
\begin{lemma}\label{l:utile}
\begin{enumerate}
\item[\rm (i)] It is possible to choose $N$ large enough in order to have the following property: if $K_{\alpha\beta}(x- y)>\epsilon$ for some $\alpha,\beta\in S(1)$ and some $(x,y)\in Q_{\bf x}\times Q_{\bf y}$, then $K_{\alpha\beta}(a- b)>\epsilon/2$ for every $(a,b)\in Q_{\bf x}\times Q_{\bf y}$. 

\item[\rm (ii)] For $\lambda>0$, and $Q\in \mathcal{Q}_\lambda$, write
\begin{equation}\label{e:los}
L(\phi_\lambda^x ; Q) := \mathcal{H}^1\left((\phi_\lambda^x)^{-1}(0)\cap Q\right).
\end{equation}
Then, for $r_\lambda$ as above there exists a finite constant $D$, independent of $\lambda$, such that
$$
\sup_{Q\in \mathcal{Q}_\lambda} \E[L(\phi_\lambda^x ; Q)^2]\leq D, \quad \lambda>0.
$$
\end{enumerate}
\end{lemma}
\noindent\begin{proof} 
Since the proof of Point (i) is essentially the same of \cite[Lemma 7.3]{NPR}, we omit it. 

Throughout all the proof of Point (ii), we set $\sigma_\lambda(x)=\sqrt{C_\lambda(x,x)}$. Moreover, we use the convention that $c_i>0$, $i=1,2,\ldots$, always denotes a constant that is {\bf independent} of $\lambda$.

The proof of Point (ii) is now divided into several steps.
Let us fix $Q\in \mathcal{Q}_\lambda$.

\medskip

{\it Step 1}. We claim that $\big|K(x,y)-2\pi+\pi\|x-y\|^2\big|\leq c_1\|x-y\|^3$ for all $x,y\in\R^2$. Indeed, fix $x\in\R^2$, and write $\widehat{K}_x(y)=K(x,y)=\E[X(x)X(y)]$.
{By definition we have} $\widehat{K}_x(x)=\E[X(x)^2]=2\pi$ and $\nabla \widehat{K}_x(x)=\E[X(x)\nabla X(x)]=0$. Thus, by a classical Taylor expansion:
\begin{equation}\label{onestar}
\widehat{K}_x(y)=2\pi+\langle ({\rm Hess} \,\widehat{K}_x)(x)(y-x),y-x\rangle+O(\|x-y\|^3),
\end{equation}
where the big $O$ is uniform with respect to $x$, because the third partial derivatives of $K_x$ are uniformly bounded  with respect to $x$.
On the other hand, using that  $K(x,y)=(2\pi) J_0(\|x-y\|)$, one can easily compute that $({\rm Hess}\,\widehat{K}_x)(x)=-\pi\,I_2$ for all $x\in\R^2$ (with $I_2$ the $2\times 2$ identity matrix). Plugging this into (\ref{onestar}), we get the announced inequality, that is,
\[
K(x,y)=\widehat{K}_x(y)={\red 2\pi } -\pi\|x-y\|^2+O(\|x-y\|^3).
\]

\medskip

{\it Step 2}. If $\lambda$ is large enough so that $\frac12\sqrt{2\pi}\leq\sigma_\lambda(x)\leq \frac32\sqrt{2\pi}$ for all $x\in Q$, we claim that 
\begin{equation}\label{doublestar}
\left|\frac{C_\lambda(x,y)}{\sigma_\lambda(x)\sigma_\lambda(y)}-1+\frac12\|x-y\|^2\right|\leq c_2\,\eta_\lambda\|x-y\|^2+c_3\|x-y\|^3
\end{equation}
for all $x,y\in  Q$.
Indeed, fix $x\in\R^2$, and write  $\widehat{C}_{x,\lambda}(y)=\frac{C_\lambda(x,y)}{\sigma_\lambda(x)\sigma_\lambda(y)}$. We have
$\widehat{C}_{x,\lambda}(x)=1$ and $\nabla \widehat{C}_{x,\lambda}(x)=0$. 
As a consequence, using in particular that the third partial derivatives of $\widehat{C}_{x,\lambda}$ are equal to that of $\frac{\widehat{K}_x}{2\pi}$ plus a remainder bounded by $O(\eta_\lambda)$, we can write that
\[
\big|\widehat{C}_{x,\lambda}(y)-1
-\langle ({\rm Hess}\, \widehat{C}_{x,\lambda})(x)(y-x),y-x\rangle\big|\leq {(c_4 + \eta_\lambda)\,\|x-y\|^3\le c_5\,\|x-y\|^3}.
\]
On the other hand, 
we have
\[
({\rm Hess}\, \widehat{C}_{x,\lambda})(x)
=-\frac12\,I_2+(({\rm Hess} \,\widehat{C}_{x,\lambda})(x)-\frac1{2\pi}({\rm Hess}\, \widehat{K}_x)(x))
\]
and the second term of the right-hand side is bounded by  $\eta_\lambda$
for any $x\in Q$.
The desired conclusion follows.

\medskip

{\it Step 3}.
According to Kac--Rice, one has
\begin{eqnarray}\label{2terms}
&&\E[L(\phi_\lambda^x;Q)^2]\\
&=&\int_{Q\times Q}\E[\|\nabla Y_\lambda(x)\|\|\nabla Y_\lambda(y)\|\big| Y_\lambda(x)=Y_\lambda(y)=0] p_{(Y_\lambda(x),Y_\lambda(y))}(0,0)dxdy.\notag
\end{eqnarray}
Using classical linear regression for Gaussian vectors, one can write
\[
\nabla Y_\lambda(x)=a_{x,y}Y_\lambda(x)+b_{x,y}Y_\lambda(y)+Z_{\lambda,x,y},
\]
with $a_{x,y}$, $b_{x,y}$ two deterministic vectors of $\R^2$ and 
 $Z_{\lambda,x,y}$ a Gaussian vector independent of $Y_\lambda(x)$ and $Y_\lambda(y)$.
As a consequence,
\begin{eqnarray*}
\E\big[\|\nabla Y_\lambda(x)\|^2\big| Y_\lambda(x)=Y_\lambda(y)=0\big]&=&\E[\|Z_{\lambda,x,y}\|^2]\\
&\leq& \E[\|\nabla Y_\lambda(x)\|^2]\leq \sup_{x\in Q}\E[\|\nabla Y_\lambda(x)\|^2].
\end{eqnarray*}
Similarly $\E\big[\|\nabla Y_\lambda(y)\|^2\big| Y_\lambda(x)=Y_\lambda(y)=0\big]\leq   \sup_{x\in Q}\E[\|\nabla Y_\lambda(x)\|^2]$ implying, by Cauchy-Schwarz, that 
\begin{equation}\label{1erterm}
\E[\|\nabla Y_\lambda(x)\|\|\nabla Y_\lambda(y)\|\big| Y_\lambda(x)=Y_\lambda(y)=0] \leq \sup_{x\in Q}\E[\|\nabla Y_\lambda(x)\|^2]
\leq c_6,
\end{equation}
the last bound being
due to the fact that, for any $x\in Q$:
\begin{eqnarray*}
\E[\|\nabla Y_\lambda(x)\|^2]&=&\E[\|\nabla X(x)\|^2]
+\partial_{(1,0)}^2(C_\lambda-K)(x,x)
+\partial_{(0,1)}^2(C_\lambda-K)(x,x))\\
&\leq& c_7+\eta_\lambda\leq c_8.
\end{eqnarray*}

\medskip

{\it Step 4}. Assume that $\lambda$ is large enough so that $\frac12\sqrt{2\pi}\leq\sigma_\lambda(x)\leq \frac32\sqrt{2\pi}$ for all $x\in Q$.
Using  (\ref{doublestar}) we have, for any $x,y\in Q$
\begin{eqnarray*}
p_{(Y_\lambda(x),Y_\lambda(y))}(0,0)
&=&\frac{1}{2\pi\sigma_\lambda(x)\sigma_\lambda(y)\,\sqrt{
1-
\frac{C_\lambda(x,y)^2}{\sigma_\lambda(x)^2\sigma_\lambda(y)^2}
}}
\leq \frac{c_9}{\sqrt{
1-
\frac{C_\lambda(x,y)^2}{\sigma_\lambda(x)^2\sigma_\lambda(y)^2}
}}\\
&\leq&\frac{c_{9}}{\sqrt{\left|\frac12\|x-y\|^2 - \big|
\frac{C_\lambda(x,y)^2}{\sigma_\lambda(x)^2\sigma_\lambda(y)^2}-1+\frac12\|x-y\|^2\big|\right|}}\\
&\leq&\frac{c_{10}}{\|x-y\|\sqrt{\big|1 -c_{11}\,\eta_\lambda-c_{12}\,\|x-y\|\big|}}.
\end{eqnarray*}
As a consequence, if $\lambda$ is large enough 
to ensure that $1 -c_{11}\,\eta_\lambda\geq \frac12$ then, for any $x,y\in Q$
satisfying $\|x-y\|\leq \frac1{4c_{12}}$, one has
\[
p_{(Y_\lambda(x),Y_\lambda(y))}(0,0)\leq \frac{c_{13}}{\|x-y\|}.
\]

\medskip

{\it Step 5}. Now, let us deal with the opposite situation where
$\|x-y\|\geq \frac1{4c_{12}}$.
For any $u,v\in[\frac12\sqrt{2\pi},\frac32\sqrt{2\pi}]$, we can write
\begin{equation}\label{elementary}
\big|2\pi-uv\big|=\left|\frac{ v^2(2\pi-u^2)+2\pi(2\pi-v^2)}{2\pi+uv}\right|
\leq \frac94|2\pi-u^2|+|2\pi-v^2|.
\end{equation}
Observe that $|\sigma_\lambda(x)^2-2\pi|\leq\eta_\lambda$ for all $x\in Q$ (this is an immediate fact, since $K(x,x)=2\pi$). Thus, we deduce from (\ref{elementary}) that
\[
\left|
\sigma_\lambda(x)\sigma_\lambda(y)-2\pi
\right|\leq \frac{13}4\eta_\lambda\quad\mbox{for all $x,y\in Q$}.
\]
This implies in turn that
\begin{eqnarray*}
&&1-
\frac{C_\lambda(x,y)^2}{\sigma_\lambda(x)^2\sigma_\lambda(y)^2}\\
&=&1-\left\{
\frac{K(x,y)}{2\pi}-\frac{K(x,y)-C_\lambda(x,y)}{\sigma_\lambda(x)\sigma_\lambda(y)}
-\frac{K(x,y)}{2\pi\sigma_\lambda(x)\sigma_\lambda(y)}\left(
\sigma_\lambda(x)\sigma_\lambda(y)-2\pi
\right)
\right\}^2\\
&\geq& 1-\left\{
\frac{|K(x,y)|}{2\pi}+\frac{|K(x,y)-C_\lambda(x,y)|}{\sigma_\lambda(x)\sigma_\lambda(y)}
+\frac{|K(x,y)|}{2\pi\sigma_\lambda(x)\sigma_\lambda(y)}\left|
\sigma_\lambda(x)\sigma_\lambda(y)-2\pi
\right|
\right\}^2\\
&\geq&1-\left\{
\frac{|K(x,y)|}{2\pi}+\frac2\pi\eta_\lambda
+\frac{2}{\pi}\left|
\sigma_\lambda(x)\sigma_\lambda(y)-2\pi
\right|
\right\}^2
\geq 1-\left\{
\frac{|K(x,y)|}{2\pi}+\frac{17}{2\pi}\eta_\lambda
\right\}^2\\
&=&1-\frac{K(x,y)^2}{4\pi^2}-\frac{17}{2\pi^2}\eta_\lambda|K(x,y)|-\frac{289}{4\pi^2}\eta_\lambda^2\\
&\geq &1-\frac{K(x,y)^2}{4\pi^2}-c_{14}\,\eta_\lambda.
\end{eqnarray*}
But $\inf_{x,y\in\R^2:\|x-y\|\geq\frac1{4c_{12}}}\big( 1-\frac{K(x,y)^2}{4\pi^2}\big):=m>0$.
Hence, for $\lambda>0$  so that $c_{14}\eta_\lambda\leq\frac{m}{2}$, one has
\[
1-
\frac{C_\lambda(x,y)^2}{\sigma_\lambda(x)^2\sigma_\lambda(y)^2}\geq c_{15}\quad\mbox{for
all $x,y\in Q$ such that $\|x-y\|\geq\frac1{4c_{12}}$},
\]
implying in turn that 
\[
p_{(Y_\lambda(x),Y_\lambda(y))}(0,0)\leq c_{16}
\] 
for all $x,y\in Q$ satisfying 
$\|x-y\|\geq\frac1{4c_{12}}$.

\medskip

{\it Step 6}. 
By merging the conclusions of steps 4 and 5, we arrive at
\[
p_{(Y_\lambda(x),Y_\lambda(y))}(0,0)\leq c_{17}\left(\frac{1}{\|x-y\|}{\bf 1}_{\{\|x-y\|\leq\frac1{4c_{12}}\}}+{\bf 1}_{\{\|x-y\|\geq\frac1{4c_{12}}\}}\right).
\]
Plugging this and (\ref{1erterm}) into (\ref{2terms}) leads to
\begin{eqnarray*}
\E[L(\phi_\lambda^x,Q)^2]&\leq& c_{17}
\int_{Q\times Q}\left(
 \frac{1}{\|x-y\|}{\bf 1}_{\{\|x-y\|\leq\frac1{4c_{12}}\}}+{\bf 1}_{\{\|x-y\|\geq\frac1{4c_{12}}\}}\right)dxdy\\
 &=&c_{17}\int_Q dx\int_{-x+Q}dz
 \left(
 \frac{1}{\|z\|}{\bf 1}_{\{\|z\|\leq\frac1{4c_{12}}\}}+{\bf 1}_{\{\|z\|\geq\frac1{4c_{12}}\}}\right)\\
 &\leq &
 c_{17}\int_Q dx\int_{B(0,\frac1{4c_{12}})}
 \frac{dz}{\|z\|}+
 c_{17}\,{\rm Leb}(Q)^2= c_{18}\,{\rm Leb}(Q)+ c_{17}{\rm Leb}(Q)^2,
\end{eqnarray*}
which automatically implies the desired conclusion.
\end{proof}

\bigskip

From now on, $N$ is fixed in such a way that the property at Point~(i) is verified. 

\begin{lemma}\label{l:comptage} As $\lambda\to \infty$, the number of singular pairs $(Q,Q')\in \mathcal{Q}_\lambda\times \mathcal{Q}_\lambda $ is $o(r_\lambda^2\log r_\lambda)$. 
\end{lemma}
\noindent\begin{proof} For every fixed $Q\in \mathcal{Q}_\lambda$, write $Z_\lambda(Q)$ for the number of cubes in $\mathcal{Q}_\lambda$ that are singular to $Q$. Then, for an absolute constant $A$ uniquely depending on $N, \epsilon$, one has that
$$
Z_\lambda(Q)\leq A\, \max_{\alpha, \beta\in S(1)} \int_{B_{2r_\lambda} - B_{2r_\lambda}} |K_{\alpha\beta}(x)|^6 dx\leq A\, \max_{\alpha, \beta\in S(1)} \int_{B_{4r_\lambda}} |K_{\alpha\beta}(x)|^6 dx.
$$
Applying \cite[Lemma 7.6]{NPR}, via an appropriate change of variables, yields that
$$
\int_{B_{4r_\lambda}} |K_{\alpha\beta}(x)|^6 dx = o (\log r_\lambda),
$$
and the desired conclusion follows immediately.
\end{proof}

\begin{remark} \label{Rem_Useful_fromSingularPairs_for_CoroKac-Rice}

	We observe that uniformly in $x, y \in B_{ r_{ \lambda} }$ one has that 
	$$
		\E \left[ \| \nabla X(x) \|| \, \| \nabla X(y)  \| \, \big| \, X(x) = 0 = X(y) \right] 
	\quad
		\text{ and }
	\quad
		p_{ ( Y_{ \lambda }( x ) , Y_{ \lambda }( y ) ) }( 0, 0)
	$$
	are finite.
	This follows from the argument leading to relation~\eqref{1erterm} 
	for the former 
	and step~5 of Point~(ii) of Lemma~\ref{l:utile} for the latter.

\end{remark}

	In the sequel we will denote by $\Sigma := \Sigma (x, y) $, for $x, y \in \R^d$ 
	the covariance matrix of the random vector $( \nabla X( x ), \nabla X (y) , X(x), X(y) )$. 
	We further define the submatrices 
	$$
	\begin{cases}
		\Sigma_{11}(x, y) &:= \Cov(   ( \nabla X( x ), \nabla X (y) ) ), \\ 
	\Sigma_{22}(x, y) &:=  \Cov( ( X(x), X(y) ) ), \\
	\Sigma_{12} (x, y) &:= \E(  ( \nabla X( x ), \nabla X (y) )  ( X(x ) , X( y ) )^t ), \\
	\Sigma_{21} (x, y) &:= \Sigma_{12}^t (x, y)
	 \end{cases}
	 $$
	 where $A^t$ denotes the transpose of the matrix.
	 For the random vector $( \nabla Y_{ \lambda} , Y_{ \lambda} )(x)$, we analogously define
	 $\Sigma^{ (  \lambda ) }$.

\begin{remark}\label{r:ru}
	 The Gaussian vector $( \nabla X( x ), \nabla X (y) )$, $x, y \in \R^d$, conditionally on the event 
	 $\{ (X(x), X ( y ) ) = ( 0, 0) \}$, is distributed as a mean zero Gaussian
	 vector with covariance matrix 
	 $$
	 	\tilde{ \Sigma } 
	= 
		\Sigma_{11} - \Sigma_{12} \Sigma_{22}^{ -1} \Sigma_{21}.
	$$

	More explicitly for the case $d=2$, $\tilde{ \Sigma } $ equals
	 $$
	 \left(\begin{array}{cccc}
	 \frac{1}{2} - ( K_{ x_1})^2/ \rho  & - K_{ x_1} K_{ x_2} / \rho  & K_{ x_1 y_1} + K K_{ x_1} K_{ y_1} / \rho  & K_{ x_1 y_2} + K K_{ x_1} K_{ y_2} / \rho  \\
	 - K_{ x_1} K_{ x_2} / \rho  & \frac{1}{2} - ( K_{ x_2})^2 / \rho   & K_{ x_2 y_1} + K K_{ x_2} K_{ y_1} / \rho  & K_{ x_2 y_2 } + K K_{ x_2} K_{ y_2} / \rho  \\
	 K_{ x_1 y_1} + K K_{ x_1} K_{ y_1} / \rho & K_{ x_2 y_1} + K K_{ x_2} K_{ y_1} / \rho & \frac{1}{2} - ( K_{ y_1 })^2 / \rho  & - K_{ y_1 } K_{ y_2 } / \rho  \\
	 K_{ x_1 y_2} + K K_{ x_1} K_{ y_2} / \rho & K_{ x_2 y_2 } + K K_{ x_2} K_{ y_2} / \rho & - K_{ y_1} K_{ y_2} / \rho & \frac{1}{2} - ( K_{ y_2 })^2 / \rho 
	 \end{array}\right)
	 $$
	 where $K_{ x_i } := \partial_{ x_i } K(x, y)$, $K_{ y_i } := \partial_{ y_i } K(x, y)$, 
	 $K_{ x_i y_j} := \partial_{ x_i } \partial_{ y_j } K (x, y)$ for $i,j=1, 2$ 
	 and $\rho := \rho(x, y) :=  \det ( \Sigma_{ 22 } ) $.
\end{remark}

	Furthermore, we introduce the following functions:

	\begin{equation*}
\begin{split}
	F_0 (x, y) 
		&:= 
		\E \left[ \| \nabla X(x) \| \, \| \nabla X(y)  \| \, \big| \, X(x) = 0 = X(y) \right] p_{ ( X(x), X(y) ) }( 0, 0); \\
	F_{ \lambda}(x, y)
		&:=
		 \E \left[ \| \nabla Y_{ \lambda }(x) \| \, \| \nabla Y_{ \lambda }(y)  \| 
		 	\, \big| \, Y_{ \lambda }( x ) = 0 = Y_{ \lambda }( y ) \right]
		 p_{ ( Y_{ \lambda }( x ) , Y_{ \lambda }( y ) ) }( 0, 0); \\
	G_{  \lambda} (x, y)
		&:=
		\E \left[ \| \nabla X(x) \| \, \| \nabla Y_{\lambda} (y)  \| \, \big| \, X(x) = 0 = Y_{\lambda} (y) \right] 
		p_{ ( X(x), Y_{\lambda} (y) ) }( 0, 0)
	\end{split}
\end{equation*}
	and also
	\begin{equation*}
\begin{split}
	H_0 (x, y) 
		&:= 
		\E \left[ \| \nabla X(x) \| \, \big| \, X(x) = 0 \right] 
		\,
		 \E \left[  \| \nabla X(y)  \| \, \big| \,  X(y) = 0 \right] 
		p_{ X(x) } (0 ) p_{ X(y)  } ( 0); \\
	H_{ \lambda}(x, y)
		&:=
		 \E \left[ \| \nabla Y_{ \lambda }(x) \| \, \big| \, Y_{ \lambda }( x ) = 0 \right] 
		 \, 
		 \E \left[ \| \nabla Y_{ \lambda }(y)  \| \, \big| \,|  Y_{ \lambda }( y ) = 0 \right]
		 p_{  Y_{ \lambda }( x ) } ( 0 ) p_{ Y_{ \lambda }( y )  }( 0); \\
	L_{  \lambda} (x, y)
		&:=
		\E \left[ \| \nabla X(x) \| \, \big| \, X(x) = 0 \right] \,  
		\E \left[ \| \nabla Y_{\lambda} (y)  \| \, \big| \,   Y_{\lambda} (y) = 0 \right] 
		p_{ X(x) } (0 ) p_{ Y_{\lambda} (y)  }(  0).
	\end{split}
\end{equation*}

	 We collect some facts in a first lemma, namely Lemma~\ref{Lem_DetCondBRW_StrictlyPos} below. 
	In a second lemma, Lemma~\ref{Lem_Bound-Diff_KacRice} below, we combine those facts 
	together with Kac--Rice formulae to obtain useful bounds on the difference between the
	nodal lengths of two Gaussian random fields.

\begin{lemma}\label{Lem_DetCondBRW_StrictlyPos}
	Assume that $\eta_{ \lambda}
	 \to 0$, as $\lambda \to \infty$.
	Then
	\begin{itemize}
	
	\item[i)] there exists a $t_0 > 0 $ such that
	\begin{equation} \label{Eq_TheoLipschitzKacRiceAwayFromDiag_EqLowerBoundCovMat}
		\inf_{ x, y : \| x - y \| > t_0  } \det(  \tilde{ \Sigma}( x, y ) )
	> 
		0
	\end{equation}
	holds;
	
	\item[ii)] for $\alpha \in \R^{ 2d } $, $A$ a positive definite matrix denote by 
	$f( \alpha,  A )$ the density of the multivariate normal distribution. 
	We have that
	$$
		| f( \alpha,   \tilde{ \Sigma}( x, y ) ) - f( \alpha,   \tilde{ \Sigma}^{ (\lambda ) }( x, y ) ) |
	\le
		P(  \alpha ) f( \alpha, \hat{ \Sigma }( x, y ) ) \eta_{ \lambda } 
	$$ 
	where $\hat{ \Sigma}( x, y ) \in \operatorname{ Conv}( \{ \Sigma( x, y ),  \tilde{ \Sigma}( x, y ) \} )$ 
	and $P(  \alpha )$ is a polynomial in $\alpha$.
	
	\end{itemize}
\end{lemma}
\noindent\begin{proof}
	\mbox{}
	\begin{itemize}
	
	\item[i)] Using Theorem~II in \cite{Ger}, it follows that each eigenvalue of $ \tilde{ \Sigma}( x, y ) $ lies
	in an open ball around some diagonal element $\tilde{ \Sigma}( i, i)$, $1 \le i \le d$,
	where the radius is given by $R_i := \sum_{j \neq i }|  \tilde{ \Sigma}( i, j ) |$.
	
	By straightforward computations we obtain the relations:
	$$
		\rho (x, y) 
	:= 
		\det( \Sigma_{22}(x, y) 
	=
		 1 - J_0^2 ( \| x - y \| ),
 	\quad 
		K_{ x_i }
	 = 
	 	-  J_1 ( \| x - y \| ) \frac{ x_i - y_i  } { \| x - y \| }
	$$
	 and $K_{ x_i} = - K_{ y_i }$, $1 \le i \le d $.
	 Furthermore $K_{ x_i y_j }$ equals 
	 $$
		\frac{ 1}{2} \left( J_0( \| x - y \| ) - J_2 ( \| x - y \| ) \right) \frac{ (x_i - y_i ) ( x_j - y_j ) }{ \| x - y \|^2 } 
		-
		J_1 ( \| x - y \| )  \frac{ ( x_i - y_i ) ( x_j - y_j ) } { \| x - y \|^3 }.
	 $$
	For every $n \ge 0$, it is well known that $J_n ( z ) \to 0$, 
	as $z \to \infty$.
	Hence, it follows that any eigenvalue of  $\tilde{ \Sigma}( x, y ) $ is strictly positive, 
	whenever $\| x - y \| $ is large enough.
	
	\item[ii) ] Recall that the formula for the density of a multivariate normal random vector is
	$$
		\frac{ 1 }{ \left( ( 2 \pi )^{ 2d } \det( A ) \right)^{ 1/2 } }
		\exp \left( 
		\frac{ - \alpha^t \operatorname{ adj } (A)  \alpha }{ 2   \det( A ) } 
		 \right)
	$$
	where $\alpha \in \R^{ 2 d }$ and $\operatorname{ adj }( \cdot )$ is the adjugate matrix, 
	i.e., the transpose of its cofactor matrix.
	Further let 
	\begin{align}
		g_{kl}(  A ) 
	&:= 
		 - \frac{  1 }{ 2  } \partial_{ a_{ k l } }  \det( A )   
		 \\
		h_{kl} ( \alpha, A )  
	&:=
		- \frac{  1 }{ 2  } 
		\left\{ \alpha^t 
			\left(
			 \det( A )  \partial_{ a_{ k l} }  \operatorname{ adj }( A) 
		 	- \operatorname{ adj }( A )	\partial_{ a_{ k l } }  \det( A )  
			\right)
		 \alpha  \right\}			
	\end{align}
	where $\partial_{ a_{ kl} } A$, $1 \le k, l \le 2d$ is the partial derivative 
	w.r.t.\ the $A (k, l )$-th element.
	Straightforward calculations show that $\partial_{a_{ kl } } f( \alpha, A )$
	equals
	\begin{equation} \label{Eq_TheoLipschitzKacRiceAwayFromDiag_EqPartDerivDens}
		f( \alpha, A ) 
			\left(
			\frac{ g_{kl}(  A ) }{ \det(  A ) } +  \frac{ h_{kl} (  \alpha, A ) }{ \det( A )^{ 2 } }
			\right)			
	\end{equation}
	By the mean value inequality, taking derivatives w.r.t.\
	the matrix entries, one has for all $\alpha \in \R^{ 2 d }$ for some $c \in \left] 0, 1 \right[$ that
	$$
		| f( \alpha,  \tilde{ \Sigma} ) - f (\alpha,  \tilde{ \Sigma}^{ ( \lambda ) } ) |
	\le
		| \nabla f( \alpha, c  \tilde{ \Sigma} + (1 - c)  \tilde{ \Sigma}^{ ( \lambda ) } ) | 
		\; \|  \tilde{ \Sigma} -  \tilde{ \Sigma}^{ ( \lambda ) } \|_{ op }.
	$$
	Since the operator norm is bounded by the Hilbert--Schmidt norm, 
	one infers 
	$\| \tilde{ \Sigma} -  \tilde{ \Sigma}^{ ( \lambda ) } \|_{ op } \le C_1  \eta_{ \lambda }$
	for some constant $C_1$ independent of $\alpha$.
	By continuity of the determinant 
	and the uniform convergence of the covariance matrices, 
	 one has for $\lambda$ large enough that
	$$
		\det(  \tilde{ \Sigma} ) /2 
	\le
		\det( c  \tilde{ \Sigma} + (1 - c)  \tilde{ \Sigma}_{ \lambda } )
	 \le
	  	2 \det( \tilde{ \Sigma} ).
	$$
	Combining this with Point~(i), we see that 
	$ \det( c  \tilde{ \Sigma} + (1 - c)  \tilde{ \Sigma}^{ ( \lambda ) } )^k$, $k= -1,- 2 $ 
	are bounded by a constant. 
	Furthermore, since Bessel functions are bounded and determinants are monomials, 
	it follows that $g_{kl}(  c  \tilde{ \Sigma} + (1 - c)  \tilde{ \Sigma}_{ \lambda } )$ is bounded by 
	a constant times $ \eta_{ \lambda } $ and
	$h_{kl}( \alpha, c  \tilde{ \Sigma} + (1 - c)  \tilde{ \Sigma}_{ \lambda } )$ 
	is bounded by $\eta_{ \lambda }$ times a polynomial in $\alpha$.
	\end{itemize}
\end{proof}

Recall the definition of $\eta_\lambda$ in (\ref{e:berardino}) and the definition of $\zeta_\lambda$ in (\ref{eta}).
\begin{lemma}\label{Lem_Bound-Diff_KacRice} 
	For every $x, y \in B_{ r_{ \lambda} }$, such that $| x - y | > t_0$, with 
	$t_0$ as found in Lemma~\ref{Lem_DetCondBRW_StrictlyPos}, it holds that
	$$
	\begin{cases}
		 | F_0 (x, y) - F_{ \lambda } (x, y) | ,  | H_0 (x, y) - H_{ \lambda } (x, y) | 
	&=
		O(  \eta_{ \lambda } ); \\
		| F_0 (x, y)   - G_{ \lambda} (x, y) |,  | H_0 (x, y) - L_{ \lambda } (x, y) | 
	&=
		O( \zeta_{ \lambda} );
	\end{cases}
	$$
	as $\lambda \to \infty$, where the constants depend only on $t_0$ and the dimension.
\end{lemma}
\noindent\begin{proof} 
	We only prove the first relation $ | F_0 (x, y) - F_{ \lambda } (x, y) | = O(  \eta_{ \lambda } ) $, 
	since the others are proven in an analogous way.

	By Remark~\ref{Rem_Useful_fromSingularPairs_for_CoroKac-Rice} and the triangle inequality, 
	it then suffices to bound
	\begin{align} \label{Eq_LemLipschOfKac-RiceEq1}
		&
		 | p_{ (X(x), X(y) ) } ( 0, 0) - p_{ ( Y_{ \lambda} ( x), Y_{ \lambda} ( y)  ) } ( 0, 0) | 
	\end{align} 
	and
	\begin{align} \label{Eq_LemLipschOfKac-RiceEq2}
		&   
		\Big | \E \left[ \| \nabla X(x) \| \, \| \nabla X(y) \|  \, \big| \, X(x) = 0 = X(y) \right]  \\
		& \quad - 
		 \E \left[ \| \nabla Y_{ \lambda }(x) \| \, \| \nabla Y_{ \lambda }(y)  \| 
		 	\, \big| \, Y_{ \lambda }( x ) = 0 = Y_{ \lambda }( y ) \right] \Big|  \nonumber
	\end{align}
	
	We will now bound \eqref{Eq_LemLipschOfKac-RiceEq1}. 
	 By continuity of the determinant and the uniform convergence of the covariance matrices, 
	 one has for $\lambda$ large enough that
	 $$
	 	\frac{1}{2} \sqrt{ \det( \Sigma_{22}(x, y) )  }  
	  \le 
	  	\sqrt{ \det( \Sigma^{ (\lambda) } _{22}(x, y) ) }  
	  \le 
	  	\frac{ 3}{2} \sqrt{ \det(  \Sigma_{22}(x, y) )  }.
	$$
	 Applying the bound: $| 1/ \sqrt{x} - 1/ \sqrt{y} | \le | y - x | ( ( \sqrt{y}+ \sqrt{x} ) \sqrt{x y} )^{-1 }$ 
	 together with Point~(i) of Lemma~\ref{Lem_DetCondBRW_StrictlyPos}, 
	 the continuity of the determinant and the boundedness of Bessel functions we find that
	 \begin{equation*}
	 	\left|  \det( \Sigma_{22}(x, y) )^{ -1 / 2  } - \det( \Sigma^{ (\lambda) } _{22}(x, y) )^{-1 / 2 } \right|
	=
		O( \eta_{ \lambda } )
	 \end{equation*}
	 as $\lambda \to \infty$.
	 Again, from the continuity of the determinant 
	 we infer that \eqref{Eq_LemLipschOfKac-RiceEq1}  is smaller than 
	 $\gamma_1  \eta_{ \lambda }$, for some constant $\gamma_1$.
	 
	 Next we turn to \eqref{Eq_LemLipschOfKac-RiceEq2}.
	By definition $\E \left[ \| \nabla X(x) \| \, \| \nabla X(y) \| \, \big| \, X(x) = 0 = X(y) \right]$ equals
	$$
	 \int_{ \R^d \times \R^d }
	 \frac{  
	   \left(  \sum_{ i=1}^d \alpha_i^2 \right)^{1/2} \left( \sum_{ i=d +1 }^{ 2d } \alpha_i^2   \right)^{1/2}
	 } {
	  \left( ( 2 \pi )^{ 2 d } \det( \tilde{ \Sigma } (x, y) ) \right)^{1/2}  }
	\exp \left( 
		\frac{ - \alpha^t \operatorname{ adj } ( \tilde{ \Sigma }(x, y) )  \alpha }{ 2  \det( \tilde{ \Sigma } (x, y) ) } 
		 \right) 
	d \alpha
	$$
	where $\operatorname{ adj }( \cdot )$ is the adjugate matrix, i.e., the transpose of its cofactor matrix.
	Now  \\
	$
		\E\left[ \| \nabla Y_{ \lambda }(x) \| \, \| \nabla Y_{ \lambda }(y)  \| 
		\, \big| \, Y_{ \lambda }( x ) = 0 = Y_{ \lambda }( y ) \right]
	$ 
	satisfies a similar identity with $ \tilde{ \Sigma }^{ ( \lambda) }(x, y )$ instead of $\tilde{ \Sigma }(x, y )$.

	From Point~(ii) of Lemma~\ref{Lem_DetCondBRW_StrictlyPos} 
	we infer \eqref{Eq_LemLipschOfKac-RiceEq2} is less than
	$$
	 \int_{ \R^d \times \R^d }
		   \left(  \sum_{ i=1}^d \alpha_i^2 \right)^{1/2} \left( \sum_{ i= d +1 }^{ 2d } \alpha_i^2   \right)^{1/2}
		P(  \alpha ) f( \alpha,  \hat{ \Sigma }( x, y ) ) \eta_{ \lambda }  
	d \alpha
	$$
	where $ \hat{ \Sigma } := c  \tilde{ \Sigma} + (1 - c)  \tilde{ \Sigma}^{ ( \lambda ) }$, for some 
	$c \in \left] 0, 1 \right[$.
	Now the matrix $ \hat{ \Sigma }$ is symmetric and positive definite. 
	Since every moment of a multivariate normal is, up to a constant, bounded by its 
	second moment, one has that 
	$$
	 \int_{ \R^d \times \R^d }
		   \left(  \sum_{ i=1}^d \alpha_i^2 \right)^{1/2} \left( \sum_{ i= d +1 }^{ 2d } \alpha_i^2   \right)^{1/2}
		 f( \alpha,   \hat{ \Sigma}(x, y) ) 
		 P( \alpha)
	d \alpha
	$$
	 is a polynomial in  $\hat{ \Sigma}(x, y)$.
	Because Bessel functions are bounded, the above integral is bounded.
	Hence \eqref{Eq_LemLipschOfKac-RiceEq2} is $O( \eta_{ \lambda})$, for $\lambda$ large enough.
	
\end{proof}

\medskip

We eventually point out the following elementary fact: if a sequence of two dimensional random vectors $\{(U_n, V_n) : n\geq 1\}$ is such that
\begin{equation}\label{easy}
\E[(U_n - V_n)^2]\to 0 \quad \mbox{and} \quad \E[U_n^2]\to 1, 
\end{equation}
then $\E[V_n^2]\to 1$.

\subsection{The proof}\label{ss:proof}
For any subset $A\subset \R^2$ and for $W = X, Y_\lambda$ ($\lambda>0$) we write
$$
L(W ; A) := \mathcal{H}^1\left(W^{-1}(0)\cap A\right), 
$$
with the simplified notation $L(W; r) := L(W ; B_r)$. The first assertion in \eqref{e:asmvpb} already appears in \eqref{e:asmv}. In view of the elementary fact recalled in (\ref{easy}), the second assertion in \eqref{e:asmvpb} will follow immediately, once we prove that, as $\lambda\to\infty$,
\begin{equation}\label{e:task}
{\bf V}(\lambda) := {\bf Var} \, [ (L(X; r_\lambda) - L(Y_\lambda; r_\lambda))] = o(r_\lambda^2\log r_\lambda). 
\end{equation}
For every $\lambda>0$ denote by $\mathcal{B}_\lambda$ the collection of all subsets of $\R^2$ having the form $Q\cap B_{r_\lambda}$, with $Q\in \mathcal{Q}_\lambda$. A pair $(B,B') \in \mathcal{B}_\lambda\times \mathcal{B}_\lambda$ is {\it singular} if the underlying pair of squares $(Q,Q')$ is. Our starting point is the obvious decomposition
$$
{\bf V}(\lambda) = {\bf V}_s(\lambda)+ {\bf V}_{ns}(\lambda),
$$
where
$$
{\bf V}_s(\lambda):= \sum_{(B,B')\in \mathcal{B}_\lambda\times \mathcal{B}_\lambda \,\, \mbox{\small singular  }}\!\!\!\!\!\!\!\!{\rm Cov}\Big\{ (L(X; B) - L(Y_\lambda; B)), (L(X; B') - L(Y_\lambda; B')) \Big\},
$$
and
$$ 
{\bf V}_{ns}(\lambda):= \sum_{(B,B') \in \mathcal{B}_\lambda\times \mathcal{B}_\lambda\,\, \mbox{\small  non singular}}{\rm Cov}\Big\{ (L(X; B) - L(Y_\lambda; B)), (L(X; B') - L(Y_\lambda; B')) \Big\}. 
$$
Since $\E[L(X,B)^2]\leq \E[L(X,Q_0)^2]$ (by stationarity) and $\E[L(Y_\lambda,B)^2]\leq D$ (by virtue of Lemma \ref{l:utile}-(i)), applying Cauchy-Schwarz and the triangle inequality yields that 
$$
{\bf V}_s(\lambda) =  O\Big( \Big| \{ (B,B')\in \mathcal{B}_\lambda\times \mathcal{B}_\lambda \, \mbox{ singular}  \} \Big| \Big),
$$
and therefore ${\bf V}_s(\lambda) = o(r_\lambda^2\log r_\lambda)$, in view of Lemma \ref{l:comptage}. We now decompose ${\bf V}_{ns}(\lambda)$ as
\begin{eqnarray*}
{\bf V}_{ns}(\lambda) &=& \sum_{(B,B')\,\, \mbox{\small  non singular}} \left[ {\rm Cov}\{ L(X; B), L(X; B')\} +{\rm Cov}\{ L(Y_\lambda; B), L(Y_\lambda; B')\}  \right. \\
&&\quad\quad\quad \quad\quad\quad \left.  - 2{\rm Cov}\{ L(X; B), L(Y_\lambda; B')\}\right] \\ 
&=:&  {\bf V}_{1}(\lambda) + {\bf V}_{2}(\lambda) - 2{\bf V}_{3}(\lambda).
\end{eqnarray*}
Now call $P_\lambda\subset B_{r_\lambda}\times B_{r_\lambda}$ the union of all cartesian product of the type $B\times B'$ such that $(B,B')$ is not singular. Exploiting the definition of non-singular pairs together with property \eqref{e:mixcov}, for $\lambda$ large enough we can represent each of the quantities  ${\bf V}_{i}(\lambda)$, $i=1,2,3$, by means of the Kac--Rice formula (or some slight variation of it --- see \cite[Chapter 6]{AW}), as follows
\begin{eqnarray*}
&&{\bf V}_{1}(\lambda) = \int_{P_\lambda} \E[\|\nabla X(x)\| \|\nabla X(y)\|\, | \, X(x) = X(y) = 0] {\red p_{ ( X(x), X(y) )}(0,0) } dxdy   \\
&& - \int_{P_\lambda} \E[\|\nabla X(x)\| \, |\, X(x)=0] \E[\|\nabla X(y)\|\, | \, X(y) = 0] {\red p_{ X(x) }(0) p_{X(y)}(0) }dxdy \\
&&=: \int_{P_\lambda}  H_1(x,y) dxdy\,;\\
&&{\bf V}_{2}(\lambda) = \int_{P_\lambda} \E[\|\nabla Y_\lambda(x)\| \|\nabla Y_\lambda(y)\|\, | \, Y_\lambda(x) = Y_\lambda(y) = 0] {\red p_{ ( Y_\lambda(x), Y_\lambda(y) ) }(0,0) } dxdy \\
&& - \int_{P_\lambda} \E[\|\nabla Y_\lambda(x)\| \, |\, Y_\lambda(x)=0] \E[\|\nabla Y_\lambda(y)\|\, | \, Y_\lambda(y) = 0] {\red p_{Y_\lambda(x)}(0) p_{Y_\lambda(y)} (0) } dxdy \\
&&=:\int_{P_\lambda}  H_2(x,y) dxdy\,;\\
&& {\bf V}_{3}(\lambda) = \int_{P_\lambda} \E[\|\nabla X(x)\| \|\nabla Y_\lambda(y)\|\, | \, X(x) = Y_\lambda(y) = 0] {\red p_{ ( X(x), Y_\lambda(y) ) }(0,0) } dxdy \\
&& - \int_{P_\lambda} \E[\|\nabla X(x)\| \, |\, X(x)=0] \E[\|\nabla Y_\lambda(y)\|\, | \, Y_\lambda(y) = 0] {\red p_{X(x)}(0) p_{Y_\lambda(y)}(0) } dxdy \\
&&=: \int_{P_\lambda}  H_3(x,y) dxdy\, ,
\end{eqnarray*}
where ${\red p_{U}( 0) }$ denotes the density of $U$ in 0, and ${\red p_{ ( U,V )}(0, 0) }$ denotes the density of $(U,V)$ in $(0,0)$. {\teal Exploiting the fact that} each $P_\lambda$ only involves non{\red -}singular pairs, we wish now to prove that, for every $1\leq j\neq k\leq 3$, and for some constant $C$ independent of $\lambda$,
\begin{equation} \label{Eq_MainTheo_Order_of_Var_of_NonSingTerms}
\int_{P_\lambda}  | H_j(x,y) -H_k(x,y)|   dxdy \leq C\, \sqrt{a(\lambda)} |P_\lambda| = O\left( \left[ \frac{r_\lambda^{17}   }{\log \lambda}\right]^{1/4} r_\lambda^4\right).
\end{equation}
{\red {\teal Once such a relation is established, the proof will be finished}, since from condition \eqref{e:21} we deduce that 
\eqref{Eq_MainTheo_Order_of_Var_of_NonSingTerms} $= o( r_{\lambda }^2 \log ( r_{ \lambda} ) )$, {and by virtue of the triangle inequality it holds that
$$
	| V_{ns} ( \lambda ) |
\le
	2 | V_1 (\lambda) - V_3 ( \lambda) | + | V_1 (\lambda) - V_2 (\lambda) | .
$$}
{\teal The validity of  \eqref{Eq_MainTheo_Order_of_Var_of_NonSingTerms} for every $j,k$ is immediately deduced from Lemma~\ref{Lem_Bound-Diff_KacRice}, as well as \eqref{e:berardino} and  \eqref{e:mixcov}, taking into account the definition \eqref{def_a_lambda}.}
}

\section{Proof of technical estimates}\label{Stechnical}

\subsection{Eigenvalues of Hilbert--Schmidt operators on expanding domains}\label{ss:hsgrow}

In the present and following section, we adopt the notation and assumptions of Section \ref{ss:statcoupling}. In particular, we fix an integer $d\geq 1$, and consider two smooth covariance kernels $K,C$ on $\R^d$, such that $K$ has bounded derivatives of every order; the symbol $B_R$ indicates the open ball with radius $R$ centred at the origin. Our proof of Theorem \ref{t:wcoupling} relies on a number of results involving Hilbert--Schmidt operators associated with the class of increasing domains $B_R$, $R\geq 1$, that we will gather below.

\smallskip

\noindent\underline{\it A simple geometric fact.} For every $A\subset \R^d$ and every $n\geq 2$, we use the notation
$$
\varepsilon_n(A):= \inf\left\{\varepsilon>0 : \mbox{$A$ can be covered by $n-1$ open balls of radius $\varepsilon$ centred in $A$}  \right\}.
$$ 
We will make use of the following elementary geometric fact 

\begin{lemma}\label{l:geo} There exist  constants $0<\beta_1<\beta_2 <\infty$ such that
$$
\frac{\beta_1}{n^{1/d}} \leq \frac{\varepsilon_n(\overline{B_R} )}{R} \leq \frac{\beta_2}{n^{1/d}}, \quad \mbox{for every } n, R\geq 1.
$$

\medskip
\end{lemma}
\noindent\begin{proof}  Consider first the case $R=1$, and denote by $Q$ the hypercube with unit side centred at the origin. It is easily seen that $\varepsilon_n(Q) \leq \beta n^{-1/d}$ for some absolute constant $\beta$, and also that $ \varepsilon_n(Q)\geq \varepsilon_n(\overline{B_1})$. To conclude the proof of the upper bound, observe that there exists an $\varepsilon$-cover of cardinality $n-1$ of $\overline{B_1}$ if and only if there exists an $(R\varepsilon)$-cover of cardinality $n-1$ of $\overline{B_R}$. The lower bound follows from the following observation: for every $\varepsilon>\varepsilon_n(\overline{B_R})$, we have that
$$
{\rm Vol}(\overline{B_R})\leq (n-1) {\rm Vol} ( B_{\varepsilon}).
$$ 

\end{proof}

\medskip

\noindent\underline{\it A class of Hilbert--Schmidt operators.} Given a finite set $J$, we denote by $L_J^2({B_R} )$ the Hilbert space of (equivalence classes of) square-integrable $\R^J$-valued functions on $B_R$, endowed with the following inner product: for $f = \{f_j : j\in J\}, g= \{g_j : j\in J\} \in  L_J^2({B_R} )$,
$$
\langle f,g \rangle_{L_J^2({B_R} )} = \sum_{j\in J} \int_{B_R} f_j(x) g_j(x) dx.
$$
Recalling the definition of $S(M)$ given in \eqref{e:sm}, for $M=0,1,...$ and $R\geq 1$, we define the operator
\begin{eqnarray}\label{e:kmr}
K^{(M,R)} &:&  L^2_{S(M)}({B_R})  \rightarrow L^2_{S(M)}({B_R})\\
&:& f = \{f_\alpha : \alpha\in S(M)\} \mapsto K^{(M,R)}f = \Big\{ (K^{(M,R)}f )_\alpha : \alpha\in S(M)\Big\},\notag
\end{eqnarray}
with
$$
(K^{(M,R)}f )_\alpha(x) := \sum_{\beta\in S(M)} \int_{B_R} K_{\alpha\beta}(x,y) f_\beta(y)\, dy, \quad x\in B_R.
$$
Note that the quantity $(K^{(M,R)}f )_\alpha(x)$ is also unambiguously defined for every $x$ in the complement of $\overline{B_R}$. It is easily checked that, for every choice of $M,R$ as above, $K^{(M,R)}$ is compact, self-adjoint, Hilbert--Schmidt and positive definite.

The eigenvalues of $K^{(M,R)}$ are denoted by 
$$
\lambda_1^{(M,R)}\geq \lambda_2^{(M,R)} \geq \cdots \geq 0.
$$
We record the following consequence of the matrix-valued Mercer's theorem (see e.g. \cite[Theorem 4.1]{dVUV}).

\begin{proposition}\label{p:mercer} Under the above notation and assumptions, let
$$
\big\{e^j =\{e^j_\alpha : \alpha\in S(M)\} : j\geq 1 \big\}
$$
be an orthonormal basis of $\overline{{\rm Im} \, K^{(M,R)}}$ (that is, of the closure of the image of $K^{(M,R)}$ in $L^2(B_R)$) composed of continuous eigenfunctions of $K^{(M,R)}$. Then, for every $\alpha, \beta \in S(M)$,
$$
K_{\alpha\beta}(x,y) = \sum_{j=1}^\infty \lambda_j^{(M,R)}e_\alpha^j (x) e_\beta^j(y),
$$
where the series converges uniformly on every compact subset of $B_R\times B_R$. This implies that $K^{(M,R)}$ is trace class and
\begin{equation}\label{e:trace1}
{\rm Tr}\, K^{(M,R)} = \sum_{j=1}^\infty \lambda_j^{(M,R)} = \sum_{\alpha\in S(M)} \int_{B_R} K_{\alpha\alpha}(x,x)\, dx.
\end{equation}
\end{proposition}

\smallskip

\noindent\underline{\it Eigenvalue decay.} Under the above introduced notation, we will now prove a useful estimate, yielding an upper bound on the decay of the eigenvalues of $\KMR$. Observe that bounds on the eigenvalue decay for kernel operators such as $\KMR$ are already available in the literature -- see e.g. the classical (and somehow definitive) reference \cite{kuhn86}, as well as \cite{hk}. However, such estimates are typically provided for kernels defined on {\it fixed domains}, whereas the applications developed in the present paper require bounds on eigenvalues where the dependence on the `increasing radius' $R$ appears explicitly. The next lemma, whose proof is partially inspired by the arguments developed in $\cite{hk}$, shows that the dependence on the radius $R$ is indeed sub-algebraic.
\begin{lemma}[\bf Eigenvalue decay for kernel operators] 
For every $M=0,1,...$ and every $\ell\geq 0$, there exists a constant $A = A(M,d,\ell)\in (0,\infty)$, exclusively depending on $M,d,\ell$, such that
\begin{equation}\label{e:eigenbound}
{ \lambda_n^{(M,R)}} \leq  A \, \frac{R^{\ell+1 +d}}{n^{\frac{\ell+1}{d}}} , \quad n\geq 1.
\end{equation}
\end{lemma} 
\noindent\begin{proof} Since $\KMR$ is a positive self-adjoint operator, one has that, {\red by Theorem~2.1 in Chapter~II of \cite{GohKrein69} } 
$$
\lambda_n^{(M,R)} = \inf\{ \| \KMR - U\|_{op} : {\rm rank} \, U<n\},
$$
where the compact notation indicates that the infimum is taken over all linear operators $U : L^2_{S(M)}({B_R})  \rightarrow L^2_{S(M)}({B_R})$ having a rank strictly less than $n$. For every smooth mapping $\psi$ on $\R^d$ every $x_0\in \R^d$ and every $\ell = 0,1,...$, we denote by 
$$
x\mapsto P_\ell(\psi, x_0)(x)
$$
the Taylor polynomial of order $\ell$ associated with $\psi$ and centered at $x_0$. Observe that $\mapsto P_\ell(\psi, x_0)$ involves a number $r=r(d,\ell)$ of monomials (each centred at $x_0$), such that the number $r$ only depends on $d$ and $\ell$. Now fix $n$ and take $\varepsilon>\varepsilon_n(\overline{B_R})$. Then, there exists a collection $\mathcal{B} = \{B(x_1, \varepsilon), ..., B(x_{n-1}, \varepsilon)\}$ of open balls centered in $\overline{B_R}$ and covering $\overline{B_R}$. Let $\phi = \{\phi_1,..., \phi_{n-1}\}$ be a partition of unity subordinated to $\mathcal{B}$. It is an elementary fact (see e.g. \cite[Theorem 2.13]{Rudin}) that one can choose $\phi$ in such a way that each $\phi_j$ is supported in the ball $B(x_j, \varepsilon)$. Define the operator
\begin{eqnarray*}
U &:&  L^2_{S(M)}({B_R})  \rightarrow L^2_{S(M)}({B_R})\\
&:& f = \{f_\alpha : \alpha\in S(M)\} \mapsto Uf = \Big\{ (Uf )_\alpha : \alpha\in S(M)\Big\}\notag
\end{eqnarray*}
by the relation
$$
(Uf )_\alpha(x) := \sum_{j=1}^{n-1}\phi_j(x) P_\ell( (\KMR f)_\alpha,x_j)(x), \quad x\in B_R,
$$
and observe that, in view of the previous considerations, ${\rm rank}\, U \leq |S(M)|\, r(d, \ell)\times (n-1) := \gamma\times (n-1)$. Also, by inverting derivation and integration, one sees that
$$
P_\ell( (\KMR f)_\alpha,x_j)(x) = \sum_{\beta\in S(M)}\int_{B_R} P_\ell (K_{\alpha\beta}(\cdot, y), x_j )(x) f_\beta(y) dy,
$$
from which we infer that, for some absolute constant $A $ and every $x\in B_R$
$$
|(\KMR f)_\alpha(x) - (Uf)_\alpha(x)| \leq A\sum_{j=1}^{n-1}\phi_j(x)\| x-x_j\|^{\ell+1} \left(\sum_{\beta\in S(M)} \int_{B_R} | f_\beta (y)| dy\right).
$$
We stress that the last inequality uses the fact that, by assumption, the derivatives of $K$ are bounded. Now, by virtue of the properties of $\phi$ one has the estimate $\sum_{j=1}^{n-1}\phi_j(x)\| x-x_j\|^{\ell+1}\leq \varepsilon^{\ell+1}$ (for every $x\in B_R$), and applying Cauchy-Schwarz we finally obtain that
$$
\| \KMR - U\|_{op}\leq A \varepsilon^{\ell+1} R^d.
$$
Letting $\varepsilon$ converge to $\varepsilon_n(\overline{B_R})$, applying Lemma \ref{l:geo} and using the fact that ${\rm rank}\, U \leq \gamma\, (n-1)$ (where the integer $\gamma$ has been defined above, and only depends on $M,d,\ell$), we see that \eqref{e:eigenbound} holds for integers of the form $n = \gamma k$, where $k\geq 1$. The fact that the conclusion is indeed true for every positive integer $n$ follows from standard arguments.

\end{proof}

\medskip

\noindent\underline{\it Square roots.} For every $M=0,1,...$ and $R\geq 1$, we define $\sqrt{\KMR}$ to be the Hilbert--Schmidt operator on $L^2_{S(M)}(B_R)$ defined by the following relations: if $\{e^j : j\geq 1\}$ denotes a continuous orthonormal basis of $\overline{{\rm Im}\, \KMR}$, then, for every $h\in  \LBR$
$$
\sqrt{\KMR} {\red h : }=  \{ (\sqrt{\KMR} h)_\alpha : \alpha\in S(M)\},
$$
with 
$$
 (\sqrt{\KMR} h)_\alpha(x) = \sum_{j=1}^\infty \sum_{\beta\in S(M)} \sqrt{\lambda_j^{(M,R)}} e_\alpha^j(x)\int_{B_R} e^j_\beta (y) h(y)\, dy,
$$
with convergence in $L^2(B_R)$. If $C$ denotes the other covariance kernel of class $\cC^{\infty, \infty}$ introduced in Section \ref{ss:statcoupling}, we define the operator $C^{(M,R)}$ (for every $M=0,1,...$ and $R\geq 1$) in the same way as above, and write
$$
\gamma^{(M,R)}_1\geq \gamma^{(M,R)}_2\geq \cdots \geq 0, 
$$
to indicate the sequence of its eigenvalues. We also define analogously the square root $\sqrt{C^{(M,R)}}$. The next estimate is crucial for our arguments.

\begin{proposition}\label{p:hssqrt} For $K,C$ as above, for every $M=0,1,...$ and $R\geq 1$, one has that
\begin{eqnarray}
\| \SKMR - \SCMR\|^2_{H.S.}&\leq & | {\rm Tr}\, \KMR - {\rm Tr}\, \CMR| + \label{e:fin} \\
&&\quad\quad\quad 2 \| \KMR - \CMR\|_{H.S.}^{1/2}\cdot {\rm Tr} \, \SKMR\notag
\end{eqnarray}
and moreover
\begin{equation}\label{e:trace2}
{\rm Tr} \, \SKMR \leq B \cdot R^{\frac{3d+1}{2}},
\end{equation}
where $B = B(M,d)$ is a finite constant uniquely depending on $M,d$.
\end{proposition}
\noindent\begin{proof} For the rest of the proof, we use the notation
$$
\alpha :=  | {\rm Tr}\, \KMR - {\rm Tr}\, \CMR|, \quad \beta:=  \| \KMR - \CMR\|_{H.S.}.
$$
Then $\beta \geq \|\KMR-\CMR\|_{op}\geq \| \SKMR-\SCMR\|^2_{op}$.  
{\red Since we have compact operators, this follows from the finite dimensional case. For a proof of the latter, see \red for instance Theorem~V.1.9. in \cite{Bha97}.} Writing moreover $D := \SKMR - \SCMR$ and choosing an orthonormal basis $\{e^j : j\geq1\}$ diagonalizing $\KMR$, we infer that
\begin{eqnarray*}
\|D\|^2_{H.S.} &=& \|\SKMR \|^2_{H.S.}+\|\SCMR \|^2_{H.S.} - 2\langle \SKMR, \SCMR\rangle_{H.S.}\\
&\leq& \alpha + 2\left|   \langle D, \SKMR\rangle_{H.S.}\right| \\
&=& \alpha+2\left| \sum_{j\geq 1} \langle De^j, \SKMR e^j\rangle_{\LBR}\rangle\right| \\
&\leq & \alpha +2 \|D\|_{op} \sum_{j\geq 1} \| \SKMR e^j\|_{\LBR}\\
&=& \alpha+2\|D\|_{op}\, {\rm Tr} \, \SKMR.
\end{eqnarray*}
The estimate \eqref{e:trace2} follows from \eqref{e:eigenbound}, by selecting $\ell = 2d$. \end{proof}

\subsection{Proof of Theorem \ref{t:wcoupling}}  \label{ss:coupling}

Fix $M\geq 0$ and $R\geq 1$, and consider a probability space $(\Omega_0, \mathscr{F}_0, \P_0)$ supporting an isonormal Gaussian process over $\LBR$, written
$$
G = \{G(h) : h\in \LBR\}.
$$
Recall that, by definition, $G$ is a centred Gaussian family, indexed by the elements of $\LBR$ and such that
$$
\E_0[G(h)G(h')] = \langle h,h'\rangle_{\LBR};
$$
observe that $G$ exists for every choice of $M$ and $R$ -- see \cite[Proposition 2.2.1]{NP}. We denote by $\{e^j : j\geq 1\}$ and $\{g^j : j\geq 1\}$ two orthonormal systems of continuous functions, composed respectively of eigenfunctions of $\KMR$ and $\CMR$ with non-zero eigenvalues; we assume that $e^j$ has eigenvalue $\lambda_j^{(M,R)}$ and $g^j$ has eigenvalue $\gamma_j^{(M,R)}$. For every $N\geq 1$ and $x\in B_R$, we set
\begin{eqnarray*}
X^N(x) &:=& \sum_{j=1}^N \sqrt{\ljmr}G(f^j ) f^j(x), \\ 
Y^N(x) &:=& \sum_{j=1}^N \sqrt{\gjmr}G(g^j ) g^j(x).
\end{eqnarray*}
Now, for every fixed $x\in B_R$ and $\alpha\in S(M)$, the sequence 
$$
X_\alpha^N(x) = \sum_{j=1}^N \sqrt{\ljmr}G(f^j ) f_\alpha^j(x), \quad N\geq 1,
$$
is Cauchy in $L^2(\P_0)$ and a.s.-$\P_0$ (by virtue of L\'evy Theorem -- see \cite[p. 23]{NP}), since $K_{\alpha\alpha}(x,x) = \sum_{j=1}^\infty \ljmr f_\alpha^j(x)^2$, were we have applied Proposition \ref{p:mercer}. Denote by $\hat{X}_\alpha(x)$ the limit of $X_\alpha^N(x)$. The following facts can be verified by a standard application of the same arguments that lead to the proof of Kolmogorov's Theorem (see e.g., \cite[Appendix A.10]{NS}):

\begin{itemize}
\item[--] the process $\Big\{ \hat{X}(x) = \{\hat{X}_\alpha(x) : \alpha\in S(M)\}, \, x\in B_R\Big\}$ admits a modification $\{ {X}(x) = \{{X}_\alpha(x) : \alpha\in S(M)\}$ such that $X_0$ (where $0$ indicates the zero multi-index) is an element of $\cC_b^\infty$ with $\P_0$-probability one and, $\P_0$-almost surely for every $x\in B_R$, $X_\alpha(x) = \partial^\alpha X_0(x)$, for every $\alpha\in S(M)$;

\item[--] for every $\alpha, \beta\in S(M)$ and for every $x,y\in B_R$
$$
\E[X_\alpha(x)X_\beta(y){]} = K_{\alpha\beta}(x,y);
$$

\item[--] for every $\alpha\in S(M)$ the real valued Gaussian field $X^N_\alpha$ converges to $X_\alpha$ uniformly on every compact subset of $B_R$. 


\end{itemize}
We define the process $Y = \{ Y_\alpha : \alpha\in S(M)\}$ in a similar way. To conclude, we just observe that, by the isometric properties of $G$,
$$
\E_0\left[\| X_0 - Y_0\|^2_{\WW^{M,2}(B_R)}\right] = \| \SKMR - \SCMR\|^2_{H.S.},
$$
and apply Theorem 2, together with the estimates
\begin{eqnarray*}
\left| {\it Tr}\, \KMR - {\it Tr}\, \CMR \right| &\leq & \sum_{\alpha\in S(M)} \int_{B_R} | K_{\alpha\alpha}(x,x) - C_{\alpha\alpha}(x,x) | dx \leq  |S(M)| |B_R|\, \eta,\\
\|\KMR - \CMR\|_{H.S.}^{1/2} &=& \left\{  \sum_{\alpha, \beta\in S(M)} \int_{B_R}\int_{B_R} (K_{\alpha\beta}(x,y) - C_{\alpha\beta}(x,y))^2 dxdy\right\}^{1/4}\\
&\leq &  |S(M)|^{1/2}|B_R|^{1/2} \sqrt{\eta}.
 \end{eqnarray*}
Relation \eqref{e:bancoupling} follows from an application of the following consequence of Sobolev's embedding theorem on open subsets of $\R^d$, such as the one stated in \cite[Theorem 2.7.2]{DD}.


\begin{lemma} For every $d,p\geq 1$, every $m> j:= \left\lfloor \frac{d}{p} \right\rfloor +1$ and every $R\geq 1$, the following estimate holds for every $u\in \cC^\infty_b(B_R)$: 
\begin{equation}\label{e:se}
\| u \|_{\cC^{m-j}(B_R)} \leq A\cdot R^{m-\frac{d}{p}} \|u\|_{\WW^{m,p}(B_R)},
\end{equation}
for a constant $A$ independent of $R$ and $u$. 
\end{lemma}
\noindent\begin{proof} The inequality is true for $R=1$ and some absolute constant $A$, by virtue of Sobolev embedding. On the other hand, since $R\geq 1$, 
$$
\|u\|_{\cC^{m-j}(B_R)} = \max_{|\alpha|\leq m-j} \| \partial^\alpha u\|_{\infty, B_R} = \max_{|\alpha|\leq m-j} \frac{1}{R^{|\alpha|}}  \| \partial^\alpha u^R\|_{\infty, B_1}\leq \|u^R \|_{\cC^{m-j}(B_1)} ,
$$
where $u^R(x) := u(Rx)$. The conclusion follows from the relation
$$
\|u^R \|_{\WW^{m,p}(B_1)}\leq R^{m- \frac{d}{p}} \| u\|_{\WW^{m,p}(B_R)},
$$
which is a consequence of the equality, valid for every $\alpha\in S(M)$
$$
\left\{ \int_{B_1} | \partial^\alpha u^R(x) | dx \right\}^{1/p} = R^{|\alpha|-\frac{d}{p}} \left\{ \int_{B_R} |\partial^\alpha u(y)  |^p dy\right\}^{1/p},
$$
deduced from the change of variables $y=Rx$.
\end{proof}

\section{Proof of Theorem \ref{t:arw}}\label{s:proofarw}

For every $n\in S$, write
$$
\widetilde{T}_n(x) := T_n\left(\frac{x}{2\pi\sqrt{n}}\right), 
$$
in such a way that $$\widetilde{u}_n(x-y) := \E\left[\widetilde{T}_n(x)\widetilde{T}_n(y)\right] =\int_{\mathbb{S}^1} e^{i\langle x-y , z \rangle} \mu_n(dz).$$ For every $n\in S$, we also denote by $\mu^\#_n$ the probability measure on the interval $[0, 2\pi]$ characterised by the following relation: for every bounded and measurable $f$,
$$
\int_{\mathbb{S}^1} f(z) \mu_n(dz) = \int_0^{2\pi} f(e^{i\theta})\mu^\#_n(d\theta);
$$
we also use the symbol $\mathbf{u}(dx)$ for the uniform {\blues probability} measure on $[0,2\pi]$. We measure the discrepancy between $\mu^\#_n$ and ${\bf u}$ in two ways: (i) {\bf Kolmogorov distance},
$$
{\bf Kol}( \mu^\#_n ,{\bf u} ) := \sup_{t\in [0,2\pi] } \left| \mu^\#_n[0,t] - \frac{t}{2\pi}\right|,
$$
and (ii) the $1$-{\bf Wasserstein distance} ${\bf W}_1$ defined as
$$
{\bf W}_1( \mu^\#_n ,{\bf u} ) = \sup_{h\in {\rm Lip}(1)} \left| \int^{2\pi}_0 h(\theta) \mu^\#_n(d\theta) - \int^{2\pi}_0 h(\theta) d\theta\right|, 
$$
where ${\rm Lip}(1)$ denotes the class of all real-valued $1$-Lipschitz functions. It is a well known fact that
$$
{\bf W}_1( \mu^\#_n ,{\bf u} ) = \int_0^{2\pi}  \left| \mu^\#_n[0,t] - \frac{t}{2\pi}\right| dt,
$$
and therefore ${\bf W}_1( \mu^\#_n ,{\bf u} )\leq 2\pi {\bf Kol}( \mu^\#_n ,{\bf u} )$.

\medskip

The conclusion of Theorem \ref{t:arw} follows from the next proposition, once we observe that, {by a change of variable},
$$
\mathcal{L}_n(\alpha_n)\stackrel{law}{=} \frac{1}{2\pi\sqrt{n}} \times {\rm length} \left( \widetilde{T}_n^{-1}(0) \cap B_{2\pi\alpha_n}\right).
$$
To simplify the discussion, for $r>0$ we also write
$$
\widetilde{ \mathcal{L}} _n(r) := {\rm length} \left( \widetilde{T}_n^{-1}(0) \cap B_{r }\right). 
$$

\begin{proposition}\label{p:arw} Fix $\rho <\frac{1}{2}\log\frac{\pi}{2}$. Then, there exists a density one sequence $\{n_j\}\subset S$ such that, as $n_j\to \infty$,
\begin{enumerate}
\item[\rm 1.] $\mathcal{N}_{n_j}\to \infty$ and $\mu_{n_j}$ converges weakly to the uniform measure on $\mathbb{S}^1$;

\item[\rm 2.] for every sequence $n\mapsto\alpha_n$ such that $\alpha_n\to \infty$ and $\alpha_n = o( (\log n)^{\rho/9})$,
$$
{\bf Var}( \widetilde{ \mathcal{L}} _n(\alpha_{n_j})) \sim \frac{\alpha_{n_j}^2\log \alpha_{n_j}}{256},
$$
and 
$$
\frac{\widetilde{ \mathcal{L}} _n(\alpha_{n_j}) - \E[\widetilde{ \mathcal{L}} _n(\alpha_{n_j})]}{\sqrt{\var(\widetilde{ \mathcal{L}} _n(\alpha_{n_j}))}}\stackrel{law}{\longrightarrow} Z\sim \mathcal{N}(0,1).
$$

\end{enumerate}
\end{proposition}
\noindent\begin{proof} Fix $\rho <\frac{1}{2}\log\frac{\pi}{2}$ and $\alpha_n = o( (\log n)^{\rho/9})$ as in the statement, and select parameters $0<\gamma<\alpha:= \frac{\log \pi}{\log 2} -1$ and $0<\kappa<\beta:= \frac{\log 2}{2}$ such that $\gamma\kappa=\rho$.  Thanks to the main result in \cite{KK}, we know that there exists a density one sequence $\{n_j\}\subset S$ such that Point 1 in the statement holds, and moreover
\begin{itemize}
\item[(a)] $\mathcal{N}_{n_j} = (\log n_j)^{\beta+\epsilon_{n_j}}$, with $\epsilon_{n_j}\to 0$, in such a way that $\mathcal{N}_{n_j}\geq  (\log n_j)^{\kappa}$ for $n_j$ sufficiently large {(see also \cite{BMW})};
\item[(b)] ${\bf Kol}( \mu^\#_{n_j} ,{\bf u} ) \leq 2 \mathcal{N}_{n_j}^{-\gamma}$ and therefore, for $n_j$ large enough,
$$
{\bf Kol}( \mu^\#_{n_j} ,{\bf u} ) \leq 2\left(\frac{1}{\log n_j}\right)^{\rho}.
$$
\end{itemize}
One key observation is that, for arbitrary multi-indices $\alpha, \beta$
$$
\left|\partial^\alpha\partial^\beta( \widetilde{u}_{n_j}(x-y) - J_0(\|x-y\|))\right| \leq \|x-y\| {\bf W}_1( \mu^\#_n ,{\bf u} ),
$$
and consequently, by the previous discussion,
$$
\left|  \partial^\alpha\partial^\beta( \widetilde{u}_{n_j}(x-y) - J_0(\|x-y\|)) \right| \leq 4\pi\|x-y\| \left(\frac{1}{\log n_j}\right)^{\rho},
$$
when $n_j$ is large enough. This yields in particular that, for some absolute constant $C$,
$$
\max_{\alpha, \beta\in S(3)}\sup_{x,y\in B_{\alpha_{n_j}}} \left|  \partial^\alpha\partial^\beta( \widetilde{u}_{n_j}(x-y) - J_0(\|x-y\|)) \right|\leq C\frac{\alpha_{n_j}}{(\log n_j)^\rho } {=:} \eta(n_j)\to 0,
$$
since $\alpha_n = o( (\log n)^\rho )$. The conclusion is obtained by reproducing verbatim the proof of Theorem \ref{t:main} given in Section \ref{s:kr}, by replacing $r_\lambda$ with $\alpha_{n_j}$, and then $\eta_\lambda$ with $\eta(n_j)$ (see \eqref{e:berardino}), as well as the right-hand side of \eqref{e:mixcov} with $$ B \, \frac{ \alpha_{n_j}} { (\log n_j)^{\rho/3}},$$ where $B$ is an absolute constant, and we have exploited \cite[Theorem 5.5]{BeMa} in the case $n=2$.
\end{proof}

\bigskip

 \noindent{\bf Acknowledgments}.  The authors wish to thank B.\ Hanin for discussing with us several aspects of his work \cite{CH} and B.\ Keeler for sharing his work \cite{Kee18} before publication. We also thank Y. Canzani, P.\ Kurlberg, L.\ Pratelli, P.\ Rigo and M.\ Sodin for fruitful discussions. 

The work of Gauthier Dierickx has been supported by the FNR grant STARS(PUL) (R-AGR-0502-10-B)  at Luxembourg University and in part by the
Collaborative Research Center ``Statistical modeling of nonlinear
dynamic processes'' (SFB 823, Teilprojekt C2) of the German Research 
Foundation (DFG). Ivan Nourdin is supported by the FNR grant APOGee (R-AGR-3585-10) at Luxembourg University. Giovanni Peccati is supported by the FNR grant FoRGES (R-AGR-3376-10) at Luxembourg University. The research of Maurizia Rossi has been supported by the Foundation Science Mathématique de Paris and (partially) by the ANR-17-CE40-0008 project {\it Unirandom}.

\small{}

\end{document}